\documentclass{amsart}
\usepackage{amsmath, amssymb}
\usepackage[utf8]{inputenc}
\usepackage[english]{babel}

\usepackage{geometry}
\usepackage{iftex}
\usepackage{amsmath, amsfonts, graphicx, color, alltt, float}
\usepackage[colorlinks = true,
linkcolor=blue,
urlcolor=blue,
citecolor=blue,
anchorcolor=blue,
backref=page]{hyperref}
\usepackage{color,graphicx}

\usepackage{url}
\usepackage{amsfonts}
\usepackage{tikz}
\usetikzlibrary{braids}
\usepackage[normalem]{ulem} 
\usetikzlibrary{decorations.pathreplacing}
\newtheorem{theorem}{Theorem}[section]

\newtheorem{definition}{Definition}[section]
\newtheorem{lemma}{Lemma}[section]
\newtheorem{proposition}{Proposition}[section]
\newtheorem{corollary}{Corollary}[section]
\newtheorem{remark}{Remark}[section]
\usepackage{graphicx}
\usetikzlibrary{arrows.meta, decorations.markings}
\usetikzlibrary{knots}
\begin{document}
	
	\title[$H(n)$-Gordian graphs of knotoids and their quotients]{$H(n)$-Gordian graphs of knotoids and their quotients: Gromov hyperbolicity}
	\author[M.\,M.~Alikbarov, A.\,Yu.~Vesnin]{M.\,M.~Alikbarov, A.\,Yu.~Vesnin}
	
	\address{Sobolev Institute of Mathematics, Siberian Branch of the Russian Academy of Sciences, Novosibirsk, Russia}
	\email{vesnin@math.nsc.ru}
	
	\address{Sobolev Institute of Mathematics, Siberian Branch of the Russian Academy of Sciences, Novosibirsk, Russia}
	\email{alikbarov70@gmail.com}
	
	\begin{abstract}
		In this paper, we study $H(n)$-moves on knotoid diagrams and the corresponding $H(n)$-Gordian graphs. We establish that the $H(n)$-Gordian graphs of knotoids are not Gromov hyperbolic. We also consider quotient graphs of the $H(n)$-Gordian graph of knotoids by the following knotoid invariants: the height of a knotoid, the span of the affine index polynomial, the $y$-span and $z$-span of the invariant $H_K(t, y, z)$, and the span of the Jones polynomial. We prove that in all these cases, the quotient graphs are Gromov hyperbolic.
	\end{abstract}
	
	
	
	\date{\today}
	
	\thanks{This work was carried out under the state assignment of the Sobolev Institute of Mathematics, Siberian Branch of the Russian Academy of Sciences (project No. FWNF-2026-0031).}
	
	\maketitle
	
	
	\section{Introduction}
	
	Originally in knot theory, the Gordian graph was defined as the graph whose vertices correspond to knots in the $3$-sphere (isotopy classes of knots), and two vertices are adjacent if the corresponding knots admit diagrams that differ by a single crossing change. The term ``Gordian'' is due to the fact that a crossing change is an unknotting operation, meaning that the diagram of any knot can be transformed into a diagram of the unknot by a finite sequence of such moves. The Gordian distance $d_G(K, K')$ between knots $K$ and $K'$ is defined as the minimum number of crossing changes required to transform a diagram of $K$ into a diagram of $K'$, where the minimum is taken over all diagrams of $K$ and $K'$. Denote by $\mathcal K$ the set of all knots in the $3$-sphere $S^3$; then $(\mathcal K, d_G)$ is a metric space~\cite{Mu}. In~\cite[Th.~C]{GG}, it was shown that the graph $(\mathcal K, d_G)$ is not Gromov hyperbolic~\cite{Gr}. In~\cite{HU}, the notion of the Gordian complex of knots was introduced. Namely, the Gordian complex of knots is the simplicial complex whose vertices correspond to knots, and a set of $n+1$ vertices $\{ K_0, K_1, \ldots, K_n \}$ forms an $n$-simplex if $d_G(K_i, K_j) = 1$ for all distinct pairs $i, j \in \{0, 1, \ldots, n\}$. Subsequently, other unknotting operations on knot diagrams were studied in knot theory, which led to the extension of the Gordian graph concept to such moves. Properties of Gordian graphs and Gordian complexes for unknotting operations on virtual knots and welded knots were investigated in~\cite{GPV, GKPV, HKOS, HNT, KPV, ZYL}.
	
	In this paper, we study $H(n)$-moves on knotoid diagrams and the corresponding $H(n)$-Gordian graphs. In the study of large graphs and complex networks, see e.g.~\cite{BRSV, KPKVB}, much attention is paid to the question of whether a graph is Gromov hyperbolic. In the context of knot theory, it was shown in~\cite{JLM} that the $H(n)$-Gordian graph of knots is not Gromov hyperbolic. We show that an analogous property holds for the $H(n)$-Gordian graphs of knotoids. However, the situation may change when turning from Gordian graphs to their quotient graphs by various knotoid invariants. Specifically, we consider quotient graphs of the $H(n)$-Gordian graph of knotoids by the following knotoid invariants: the height of a knotoid, the span of the affine index polynomial, the $y$-span and $z$-span of the invariant $H_K(t, y, z)$, and the span of the Jones polynomial, and we prove that in all these cases the quotient graphs are Gromov hyperbolic. 
	
	The paper is organized as follows. In Section~\ref{sec2}, we recall the definition of Gromov hyperbolicity for metric spaces and the notion of the $H(n)$-Gordian graph of knots, and review known results on the structure of this graph. In Section~\ref{sec3}, we recall the concept of a knotoid introduced by Turaev in 2012, define the $H(n)$-Gordian graphs of knotoids, and show that these graphs are not Gromov hyperbolic, see Theorem~\ref{theorem1.2}. In Section~\ref{sec4}, we recall the definitions and properties of several known knotoid invariants: the height of a knotoid~\cite{Tu}, the span of the affine index polynomial~\cite{GuKa}, and the $y$-span and $z$-span of the invariant $H_K(t, y, z)$~\cite{FLV}. In Section~\ref{sec5}, we introduce an infinite family of spiral knotoids $S_m$ with $2m$ crossings in a diagram and compute the values of the invariants from Section~\ref{sec4} for them. In Sections~\ref{sec6} and~\ref{sec7}, we establish the Gromov hyperbolicity of the quotient graphs of the $H(n)$-Gordian graph of knotoids by the invariants from Section~\ref{sec4}, see Theorem~\ref{theorem4.1}, as well as by the span of the Jones polynomial of a knotoid~\cite{Tu}, see Theorem~\ref{theorem1.6}.

\section{Gromov hyperbolicity and Gordian graphs of knots} \label{sec2}

We will use the definition of hyperbolicity via the four-point condition, which was given by Gromov in~\cite{Gr}.

\begin{definition} \cite[Ch.~III.H.1, p.410]{BH} \label{def2.4}
	{\rm
		A metric space $(X,d)$ is called \textit{$\delta$-hyperbolic} if for every four points $x, y, z, w \in X$, the following inequality holds:
		\begin{equation} \label{eqn0}
			d(x,y)+d(z,w) \leq \max\{d(x,z)+d(y,w), d(x,w)+d(y,z)\} + \delta.
		\end{equation}
		A metric space is called \textit{Gromov hyperbolic} if there exists $\delta \geq 0$ such that the space is $\delta$-hyperbolic. 
	}
\end{definition}

In the study of Gromov hyperbolicity, the notion of quasi-isometry plays an essential role.

\begin{definition} \label{def2.5}
	{\rm
		Let $(X, d_X)$ and $(Y, d_Y)$ be metric spaces. A map $f:X\rightarrow Y$ (not necessarily continuous) is called a \textit{quasi-isometry} if there exist constants $A \geq 1$, $B \geq 0$, and $C \geq 0$ such that the following properties hold:
		\begin{itemize}
			\item[1.] For any two points $x, x' \in X$,
			$$
			\frac{1}{A} \cdot d_X(x,x') - B \leq d_Y(f(x),f(x'))\leq A \cdot d_X(x,x') + B.
			$$
			\item[2.] For every point $y \in Y$, there exists a point $x \in X$ such that $d_Y(f(x),y)\leq C$.
		\end{itemize}
	}
\end{definition}

\begin{proposition} \cite[Pr.~1]{JLM} \label{prop2.1}
	Let $f: (X, d_X) \rightarrow (Y, d_Y)$ be a quasi-isometry of geodesic metric spaces. If $(X, d_X)$ is Gromov hyperbolic, then $(Y, d_Y)$ is also Gromov hyperbolic.
\end{proposition}

\begin{figure}[h]
	\centering
	\begin{tikzpicture}[scale=0.9, baseline=(current bounding box.center)]
		\begin{scope}
			\draw[very thick, black, dashed] (0,0) circle (2.0);
			\foreach \a/\label in {0/1, 30/2, 60/3, 90/4, 120/5, 150/6, 240/{2n-3}, 270/{2n-2}, 300/{2n-1}, 330/{2n}} 
			{
				\node[black, circle, fill, inner sep=1.2] (L\a) at (\a:2.) {};
				\node[black] at (\a:2.5) {\small $\label$}; 
			}
			\node[black] at (185:1.6cm) {$\vdots$};
			\draw[very thick, black] (L0) to[bend left=55] (L30);
			\draw[very thick, black] (L60) to[bend left=55] (L90);
			\draw[very thick, black] (L120) to[bend left=55] (L150);
			\draw[very thick, black] (L240) to[bend left=55] (L270);
			\draw[very thick, black] (L300) to[bend left=55] (L330);
		\end{scope}
		\draw[<->, thick] (3.,0) -- (4.,0) node[midway, above] {$H(n)$};
		\begin{scope}[xshift=7.cm]
			\draw[very thick, black, dashed] (0,0) circle (2.);
			\foreach \a/\label in {0/1, 30/2, 60/3, 90/4, 120/5, 150/6, 240/{2n-3}, 270/{2n-2}, 300/{2n-1}, 330/{2n}} 
			{
				\node[circle, black, fill, inner sep=1.2] (R\a) at (\a:2cm) {};
				\node[black] at (\a:2.5cm) {\small $\label$};
			}
			\node[black] at (185:1.5cm) {$\vdots$};
			\draw[very thick, black] (R330) to[bend left=55] (R0);
			\draw[very thick, black] (R30) to[bend left=55] (R60);
			\draw[very thick, black] (R90) to[bend left=45] (R120);
			\draw[very thick, black] (R270) to[bend left=45] (R300);
			\draw[very thick, black] (R150) to[bend left=30] ++(-40:0.5cm);
			\draw[very thick, black] (R240) to[bend right=30] ++(70:0.5cm);
		\end{scope}
	\end{tikzpicture}
	\caption{An $H(n)$-move on diagrams.}  \label{fig1}
\end{figure}
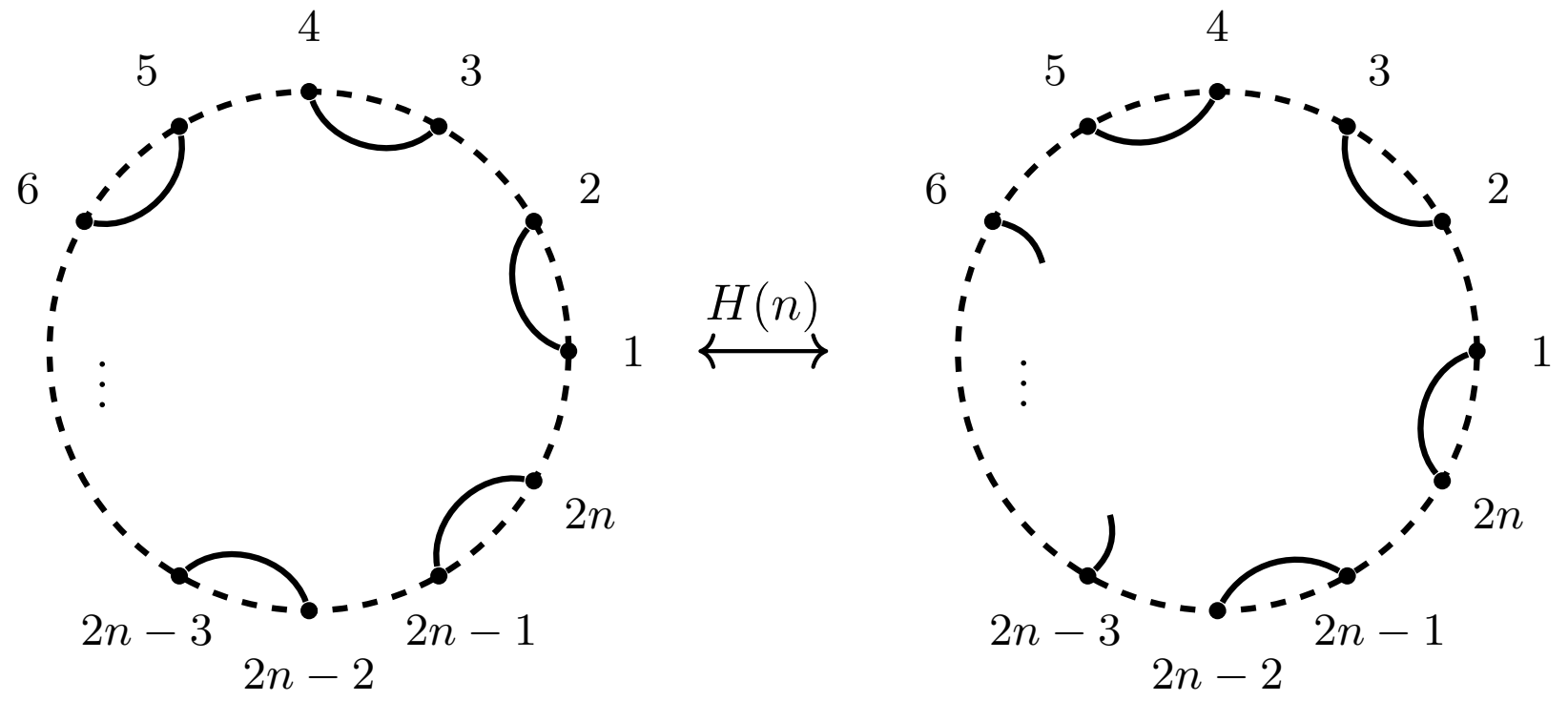

In~\cite{HNT}, for $n \geq 2$, $H(n)$-moves on knot diagrams were defined, see Fig.~\ref{fig1}, and it was proved that they are unknotting operations, meaning that a diagram of any knot can be transformed into a diagram of the unknot by a finite sequence of $H(n)$-moves, Reidemeister moves, and planar isotopies.
Here, we consider only those $H(n)$-moves that transform a knot diagram into a knot diagram. Since these moves are unknotting operations, for any two knots their diagrams can be transformed into each other by a finite number of $H(n)$-moves, Reidemeister moves, and planar isotopies. We define the $H(n)$-distance $d_{H(n)} (K, K')$ between knots $K$ and $K'$ as the minimum number of $H(n)$-moves required to transform a diagram of $K$ into a diagram of $K'$, where the minimum is taken over all diagrams of these knots. In~\cite[Lemma~1]{HNT}, it was noted that for any $n \geq 2$, an $H(n)$-move can be realized by an $H(n+1)$-move.

\begin{definition} \label{def2.2}
	{\rm
		Let $n \geq 2$. The \textit{$H(n)$-Gordian graph of knots}, denoted by $\Gamma(H(n))$, is the simple graph whose vertices are knots in $S^3$, where two knots $K$ and $K'$ are connected by an edge if some diagram of $K$ can be transformed into a diagram of $K'$ by a single $H(n)$-move, together with Reidemeister moves and planar isotopy if necessary.
	}
\end{definition}

Since an $H(n)$-move is an unknotting operation, the graph $\Gamma(H(n))$ is connected. Equipped with the metric $d_{H(n)} (K, K')$, it becomes a geodesic metric space. In this metric, the distance between two vertices coincides with the length of a shortest path connecting them in the graph. For the $H(n)$-Gordian graphs of knots, the following property holds. 

\begin{proposition} \cite[Cor.~2]{JLM} \label{prop2.2}
	For any $n \geq 2$, the $H(n)$-Gordian graph of knots $\Gamma(H(n))$ is quasi-isometric to $\Gamma(H(2))$.
\end{proposition}

In~\cite{ZYL}, it was shown that for any $0$-simplex $p$ in the $H(n)$-Gordian complex $\mathcal G_{H(n)}$, there exists a simplex $\sigma$ of arbitrarily large dimension such that $p$ is a subcomplex of $\sigma$ for all $n \geq 2$.  In~\cite{JLM}, the question of whether the $H(n)$-Gordian graph of knots is Gromov hyperbolic was investigated; namely, the following result was established.

\begin{theorem} \cite[Th.~1]{JLM}  \label{theorem1.1}
	For $n \geq 2$, the $H(n)$-Gordian graph $\Gamma(H(n))$ of knots is not Gromov hyperbolic.
\end{theorem}

Along with Gordian graphs of knots, graphs obtained from them by taking quotients by certain knot invariants are also of interest. 

\begin{definition} \label{def2.3} {\rm 
		Let $\Gamma(H(n))$ be the $H(n)$-Gordian graph of knots, and let $\operatorname{Inv}$ be a knot invariant. The \textit{quotient graph} of $\Gamma(H(n))$ by the invariant $\operatorname{Inv}$, denoted by $\Gamma (H(n), \operatorname{Inv})$, is the simple graph whose vertices are the equivalence classes of knots $[K]$, where two knots $K_1$ and $K_2$ are equivalent if $\operatorname{Inv}(K_1) = \operatorname{Inv}(K_2)$. Two classes $[K_1]$ and $[K_2]$ are connected by an edge if there exist representatives of these classes that can be transformed into each other by a single $H(n)$-move.}
\end{definition}

In the next proposition, the Alexander polynomial $\Delta(t)$ (see the definition in~\cite{Li}) is used as the knot invariant $\operatorname{Inv}$.

\begin{proposition} \label{prop1.1}
	For $n \geq 2$, denote by $\Gamma(H(n), \Delta)$ the quotient graph of the $H(n)$-Gordian graph $\Gamma(H(n))$ of knots by the Alexander polynomial $\Delta(t)$. Then the graph $\Gamma(H(n), \Delta)$ is
	\begin{enumerate}
		\item $6$-hyperbolic for $n=2$,
		\item $2$-hyperbolic for $n=3$ and $n=4$,
		\item $0$-hyperbolic for $n\geq 5$.
	\end{enumerate}
\end{proposition}

\begin{proof} 
	Denote by $[K]_{\Delta}$ the equivalence class of knots whose Alexander polynomial coincides with the Alexander polynomial of the knot $K$, and by $d_{H(n)} ([K], [K'])$ the distance in $\Gamma(H(n), \Delta)$. 
	As shown in~\cite[Th.~2]{HNT}, a crossing change can be realized by two $H(2)$-moves. In~\cite{Ko}, it was proved that for any knot $K$, there exists a knot $K' \in [K]_{\Delta}$ at a distance of one crossing change from the unknot $U$, which implies $d_{H(2)} ([K]_{\Delta}, [U]_{\Delta}) \leq 2$. Consequently, for arbitrary knots $K$ and $L$, we obtain $d_{H(2)} ([K]_{\Delta}, [L]_{\Delta}) \leq d_{H(2)} ([K]_{\Delta}, [U]_{\Delta})  + d_{H(2)} ([U]_{\Delta}, [L]_{\Delta}) \leq 4$. 
	
	In~\cite[Th.~5]{JLM}, it was shown that $m$ applications of $H(n)$-moves can be realized by a single $H((n-1)m+1)$-move. Therefore, for arbitrary knots $K$ and $L$, the following estimates hold:
	\begin{equation} \label{eqn:diameter}
		d_{H(3)} ([K]_{\Delta}, [L]_{\Delta}) \leq 2, \qquad  d_{H(4)} ([K]_{\Delta}, [L]_{\Delta}) \leq 2 \qquad  
		\text{and}  \qquad d_{H(n)} ([K]_{\Delta}, [L]_{\Delta}) \leq 1
	\end{equation} 
	for $n \geq 5$.
	
	Consider four vertices $x,y,z,w$ in $\Gamma(H(n), \Delta)$ and set
	$$
	M = \max\{d(x,z)+d(y,w), d(x,w)+d(y,z)\}.
	$$ 
	
	If $M = 0$, then $x=z$, $y=w$, $w=x$, and $z=y$. In this case, the left-hand side of inequality~(\ref{eqn0}) vanishes, so inequality~(\ref{eqn0}) holds for any $\delta \geq 0$. 
	
	If $M = 1$, then without loss of generality we may assume that $d(x,z) = 1$ and $d(y,w)=0$. Since $d(x,w)+d(y,z) \leq 1$, the following cases are possible: 
	\begin{itemize}
		\item[(1)] $d(x,w) =0$ and $d(y,z)=0$; then $x=w$ and $y=z$, which together with $y=w$ gives $x=z$, contradicting $d(x,z) = 1$; hence this case is impossible. 
		\item[(2)] $d(x,w) =1$ and $d(y,z)=0$; then $y=z$, and together with $y=w$ we get $y=z=w$; hence the left-hand side of inequality~(\ref{eqn0}) is equal to $1$, and inequality~(\ref{eqn0}) holds for any $\delta \geq 0$; 
		\item[(3)] $d(x,w) =0$ and $d(y,z)=1$; then $x=w$, and together with $y=w$ we get $x=y=w$; hence the left-hand side of inequality~(\ref{eqn0}) is equal to $1$, and inequality~(\ref{eqn0}) holds for any $\delta \geq 0$.
	\end{itemize}
	
Now let $M \geq 2$. Let $D$ be the diameter of the graph $\Gamma(H(n), \Delta)$. Then the left-hand side of inequality~(\ref{eqn0}) does not exceed $2D$, while the right-hand side is at least $2+\delta$, so inequality~(\ref{eqn0}) holds for any $\delta \geq 2D-2$. By~(\ref{eqn:diameter}), the diameter $D$ is finite. Specifically, for $\Gamma(H(2), \Delta)$ we have $D \leq 4$, and to satisfy $\delta \geq 2D-2$ we may set $\delta=6$; for $\Gamma(H(3), \Delta)$ and $\Gamma(H(4), \Delta)$ we have $D \leq 2$, so we may set $\delta=2$; for $\Gamma(H(n), \Delta)$ with $n \geq 5$ we have $D \leq 1$, so we may set $\delta=0$. The proof is complete.
\end{proof}

In~\cite{FMS}, properties of the $H(2)$-Gordian graph of knots quotiented by various invariants associated with the Jones polynomial $V(t)$ are studied. Denote by $\mathbf{span}(V)$ the span of the Jones polynomial.

\begin{theorem} \cite[Th.~4.7]{FMS}  \label{theorem1.3}
	Denote by $\Gamma(H(n), \mathbf{span}(V))$ the quotient graph of the $H(n)$-Gordian graph of knots by the span $\mathbf{span}(V)$ of the Jones polynomial $V(t)$. Then the graph $\Gamma(H(n), \mathbf{span}(V))$ is $\delta$-hyperbolic for some $\delta \geq 0$.
\end{theorem}

\section{Gordian graphs of knotoids} \label{sec3}

Recall that knotoids were introduced by V.\,G.~Turaev in 2012 in~\cite{Tu}. Namely, a knotoid in an oriented surface $S$ is defined as an equivalence class of diagrams of immersions of the oriented interval $[0,1]$ into $S$ endowed with over/under-crossing information at the double points. Equivalence of diagrams is considered up to Reidemeister moves. Knotoids can be viewed as a natural generalization of knots. Many operations on knots, in particular $H(n)$-moves, carry over naturally to the class of knotoids in the $2$-sphere, also called \textit{spherical knotoids}. Henceforth, all knotoids considered are assumed to be spherical.

The definition of $H(n)$-moves for knotoid diagrams is analogous to that for knot diagrams. Here, we consider only those $H(n)$-moves that transform a knotoid diagram into a knotoid diagram. An $H(n)$-move is also an unknotting operation for knotoids. The proof for knots in~\cite{HNT} carries over verbatim. Recall that this proof implies the inequality $d_{H(n)} (K,K')\geq d_{H(n+1)}(K,K')$. Analogously to the graph $\Gamma(H(n))$, we define the $H(n)$-Gordian graph of knotoids $G(H(n))$. To distinguish the distances in $\Gamma(H(n))$ and $G(H(n))$, we will use the notation $d^{\Gamma}_{H(n)}$ in the former case and $d^G_{H(n)}$ in the latter.

Let $D \subset S^2$ be a diagram of a spherical knotoid $K$. A \emph{shortcut} for  $D$ is a path $\alpha$ in $S^2$ connecting the tail and the head that intersects $D$ transversally and does not intersect the crossings. (The endpoints are not considered to be intersections.) Following~\cite{Tu}, we associate to a diagram $D$ of $K$ two knots $K_+$ and $K_-$, whose diagrams are obtained by connecting the endpoints of $D$ with a shortcut $\alpha$, such that in the first case $\alpha$ passes over $D$ at each intersection of $\alpha$ with $D$, and in the second case $\alpha$ passes under $D$. As shown in~\cite{Tu}, the knots $K_+$ and $K_-$ do not depend on the choice of the diagram $D$ of $K$ or the choice of the shortcut $\alpha$. Let $\phi$ be the map that associates to a knot $K$ the knotoid $\phi(K)$ obtained by deleting a small neighborhood of a point on a knot diagram that is not a crossing. In~\cite{Tu}, it was shown that $\phi$ is injective. We say that $\phi(K)$ is a \textit{knot-type knotoid}. Clearly, $\phi(K)_+ = \phi(K)_- = K$, since one can choose as $\alpha$ the arc removed under the map $\phi$, which by definition does not intersect the diagram of $\phi(K)$. 

\begin{lemma} \label{lemma3.1}
	Let $K$ and $K'$ be two knotoids in the sphere $S^2$. Then for any $n \geq 2$, the following inequality holds:
	$$
	d^G_{H(n)}(K,K') \geq \max \{d^{\Gamma}_{H(n)}(K_-, K'_-), d^{\Gamma}_{H(n)}(K _+, K'_+)\}.
	$$
\end{lemma}

\begin{proof}
	Suppose diagrams $D$ and $D'$ of the knotoids $K$ and $K'$ are related by a sequence of $H(n)$-moves. Construct diagrams of the knots $K_-$ and $K'_-$. Apply $H(n)$-moves to the diagram of $K_-$ in the same regions as for $D$. If the arc $\alpha$ connecting the endpoints of $D$ enters the region of a move, we push it outside the region using the second and third Reidemeister moves. This is always possible because $\alpha$ passes under $D$. Thus, transforming the diagram of $K_-$ into the diagram of $K'_-$ requires no more $H(n)$-moves than transforming $D$ into $D'$. Consequently, the distance between the knots $K_-$ and $K'_-$ in $\Gamma(H(n))$ does not exceed the distance between the knotoids $K$ and $K'$ in $G(H(n))$. The case of the knots $K_+$ and $K'_+$ is treated analogously.
\end{proof}

\begin{lemma} \label{lemma3.2}
	Let $K$ and $K'$ be two knots in $S^3$, and let $\phi$ be the map sending a knot to a knot-type knotoid. Then $d^{\Gamma}_{H(n)}(K,K') = d^G_{H(n)}(\phi(K),\phi(K'))$.
\end{lemma}

\begin{proof}
	Suppose the knots $K$ and $K'$ correspond to distinct vertices of the graph $\Gamma(H(n))$ and $d^{\Gamma}_{H(n)}(K,K') = s > 0$. We show that a shortest path $K = K_0, \ldots, K_s = K'$ connecting these vertices in $\Gamma(H(n))$ maps to a shortest path $\phi(K) = \phi(K_0), \ldots, \phi(K_s) = \phi(K')$ connecting the vertices in $G(H(n))$ corresponding to the knotoids $\phi(K)$ and $\phi(K')$.
	
	First, let us verify that the path $\phi(K_0), \ldots, \phi(K_s)$ lies in $G(H(n))$. Indeed, since the map $\phi$ does not depend on the choice of the point whose neighborhood is removed, in each $K_i$ we can choose such a point so that some neighborhood of it lies outside the region where the $H(n)$-move is applied.
	
	Now let us verify that this path is shortest. Suppose to the contrary that in $G(H(n))$ there exists a path $N_{0}, \ldots, N_t$ connecting the knotoids $N_0 = \phi(K)$ and $N_t = \phi(K')$ such that $t < s$. Then $d^G_{H(n)}(\phi(K),\phi(K')) < s = d^{\Gamma}_{H(n)}(K,K')$. From the definition of the map $\phi$, it follows that $\phi(K)_-$ and $\phi(K')_-$ are the knots $K$ and $K'$, respectively. Therefore, by Lemma~\ref{lemma3.1}, we have
	$$
	t \geq d^G_{H(n)}(\phi(K),\phi(K')) \geq d^{\Gamma}_{H(n)}(K,K') = s.
	$$ 
	This contradicts the assumption $t < s$. Thus, $t = s$, and the path $\phi(K_0), \ldots, \phi(K_s)$ has the minimum length among all paths in $G(H(n))$ connecting the knotoids $\phi(K)$ and $\phi(K')$.
\end{proof}

\begin{theorem}\label{theorem1.2}
	Let $G(H(n))$ be the $H(n)$-Gordian graph of knotoids, where $n \geq 2$. Then the graph $G(H(n))$ is not Gromov hyperbolic.
\end{theorem}

\begin{proof} 
	We show that for any $\delta \geq 0$, there exists a quadruple of vertices $x, y, z, w$ in $G(H(n))$ for which inequality~(\ref{eqn0}) fails. We look for such a quadruple among the vertices corresponding to knot-type knotoids. Since by Theorem~\ref{theorem1.1} the graph $\Gamma(H(n))$ is not Gromov hyperbolic, for any $\delta \geq 0$ there exists a quadruple of vertices $a$, $b$, $c$, $d$ in $\Gamma(H(n))$, corresponding to some knots $K_1$, $K_2$, $K_3$, $K_4$, such that the inequality
	\begin{equation} \label{eqn102}
		d(a,b)+d(c,d) \leq \max\{d(a,c)+d(b,d), d(a,d)+d(b,c)\} + \delta  
	\end{equation}		
	does not hold. Associate to these knots the knot-type knotoids $\phi(K_1)$, $\phi(K_2)$, $\phi(K_3)$, $\phi(K_4)$. Choose $x, y, z, w \in G(H(n))$ to be the vertices corresponding to these knotoids. Since by Lemma~\ref{lemma3.2} we have $d_{H(n)}^G (\phi(K), \phi(K')) = d_{H(n)}^{\Gamma} (K, K')$, the required property holds.
\end{proof}

\section{Some invariants of knotoids} \label{sec4}

\textbf{4.1. Height of a knotoid.} Let $\alpha$ be a shortcut for a knotoid diagram $D \subset S^2$. Of course, every knotoid diagram has many shortcuts. For each diagram $D$ and shortcut $\alpha$ the \emph{height} $h(D,\alpha)$ is the number of intersection points with $\alpha$.

\begin{definition} \label{def3.1}{\rm 
For a spherical knotoid $K$ the \emph{height} $h(K)$ is the minimum of the heights $h(D, \alpha)$ over all its diagrams $D$ and shortcuts $\alpha$.
}
\end{definition}

It follows from the definition that the height is an invariant of a spherical knotoid. The height of knotoid was called \textit{complexity} in~\cite{Tu}.


Following~\cite{Tu}, we recall the definition of the connected sum of spherical knotoids.

\begin{definition} \label{def3.2}{\rm 
		Consider diagrams $D$ and $D'$ of spherical knotoids $K$ and $K'$ (each in its own sphere). Near the head of the diagram $D$, cut out a small disk neighborhood whose boundary intersects $D$ in a single point. Similarly, cut out a small disk neighborhood near the tail of the diagram $D'$. Glue the resulting punctured spheres along their boundaries, identifying the points of $D$ and $D'$ lying on the boundaries. We obtain a new diagram $D\# D'$, whose equivalence class is called the \textit{connected sum} $K \# K'$ of the knotoids $K$ and $K'$.}
\end{definition}

Clearly, the connected sum $K \# K'$ is well-defined and does not depend on the choice of diagrams $D$ and $D'$. Moreover, a spherical knotoid $K$ is unchanged by taking the connected sum with the trivial knotoid $U$: $K\#U = K$ and $U \# K = K$. In~\cite{Tu}, the following property was established.

\begin{theorem} \cite[Th.~4.3]{Tu} \label{theorem3.1}
	For any spherical knotoids $K$ and $K'$, the equality $h(K\#K') = h(K)+h(K')$ holds, where $K\#K'$ is the connected sum of $K$ and $K'$.
\end{theorem}

\textbf{4.2. Affine index polynomial.} The affine index polynomial $P_K(t)$, which is an invariant of a knotoid, was defined in~\cite{GuKa}.

\begin{figure}[htbp]
\begin{center}
	\begin{tikzpicture}
		\draw[black, very thick, ->] (0.4,0.6) -- (0,1);
		\node at (2.5,0.5) {$\operatorname{sgn} (c) = +1$};
		\draw[black, very thick] (0.6,0.4) -- (1,0);
		\draw[black, very thick, ->] (0,0) -- (1,1);
	\end{tikzpicture}
	\quad \quad
	\begin{tikzpicture}
		\draw[black, very thick, ->] (1,0) -- (0,1);
		\node at (2.5,0.5) {$\operatorname{sgn} (c) = -1$};
		\draw[black, very thick] (0,0) -- (0.4,0.4);
		\draw[black, very thick, ->] (0.6,0.6) -- (1,1);
	\end{tikzpicture}
\end{center}
\caption{The sign of a crossing.} \label{fig4}
\end{figure}
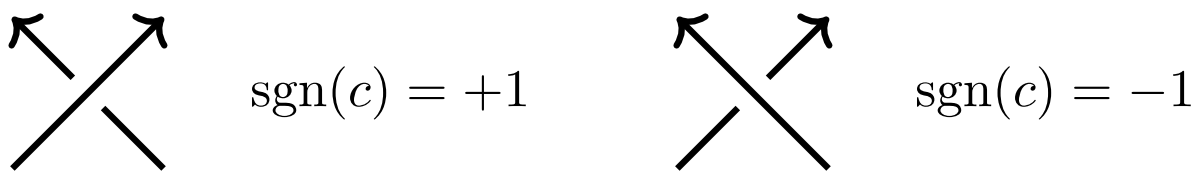

Let $c$ be a crossing of a diagram $D$ of a knotoid $K$. The function $\operatorname{sgn}(c)$ assigns to the crossing the number $+1$ or $-1$ depending on its type, as shown in Fig.~\ref{fig4}.

\begin{figure}[htbp]
\begin{center}
	\begin{tikzpicture}
		\draw[black, very thick, ->] (1,0) -- (0,1);
		\node at (-0.2,0) {$a$};
		\node at (1.2,0) {$b$};
		\node at (-0.5,1) {$b+1$};
		\node at (1.5,1) {$a-1$};
		\draw[black, very thick, ->] (0,0) -- (1,1);
	\end{tikzpicture}
\end{center}
\caption{Cheng coloring in a neighborhood of a crossing.} \label{fig5}
\end{figure}
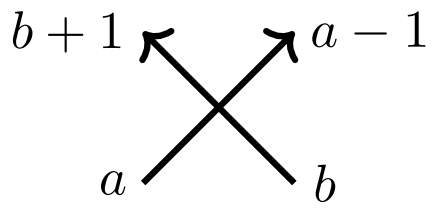

The polynomial $P_D (t)$ is constructed by assigning integer labels to the arcs of the diagram satisfying the rule shown in Fig.~\ref{fig5} in a neighborhood of each crossing. The labeling rule does not depend on the crossing type. Such a labeling is often called a \textit{Cheng coloring}.

\begin{definition} \label{def3.4} {\rm 
	Suppose a Cheng coloring is given on a diagram $D$ of a knotoid $K$. The \textit{weight} $w_D(c)$ of a crossing $c \in D$ is defined as}
$$
w_D(c) = \operatorname{sgn}(c)\cdot(a-(b+1)).
$$
\end{definition}

\begin{definition} \label{def3.5} {\rm 
	Let $D$ be a diagram of a spherical knotoid $K$. The \textit{affine index polynomial} of the diagram $D$ is defined as the polynomial
	$$
	P_D(t) = \sum_{c \in D} \operatorname{sgn}(c)\cdot(t^{w_D(c)}-1),
	$$
	where the sum is taken over all crossings.}
	\end{definition}
	
	In~\cite[Th.~4.8]{GuKa}, it was established that the affine index polynomial $P_D(t)$ is an invariant of a spherical knotoid $K$. Therefore, we henceforth denote it by $P_K(t)$. The following property was also proved there.
	
	\begin{theorem} \cite[Th.~4.12]{GuKa} \label{theorem3.2}
Let $K$ be a spherical knotoid, $h(K)$ its height, and $\operatorname{maxdeg}(P_K(t))$ the maximum degree of the polynomial $P_K(t)$. Then $h(K) \geq \operatorname{maxdeg}(P_K(t))$.
\end{theorem}

\begin{definition} \label{def3.10} {\rm 
The \textit{span} $\mathbf{span}(P_K(t))$ of the nontrivial affine index polynomial $P_K(t)$ is defined as the difference between the maximum and minimum degrees of $t$ in $P_K(t)$. For trivial polynomial we assume span equals $0$. }
\end{definition}

In~\cite{GuKa}, the following property was established.

\begin{theorem} \cite[Th.~4.10]{GuKa} \label{theorem3.3}
If $K$ is a spherical knotoid, then $P_K(t) = P_K(t^{-1})$.
\end{theorem}

Thus, if the maximum degree of the polynomial $P_K(t)$ is $n$, then $(-n)$ is its minimum degree.

\begin{corollary} \label{lemma3.4}
The span $\mathbf{span}(P_K(t))$ of the affine index polynomial $P_K(t)$ of a spherical knotoid $K$ is an even integer.
\end{corollary}

\textbf{4.3. The invariant $H_K(t,y,z)$.} In~\cite{FLV}, a transcendental function $H_K(t,y,z)$ was introduced, which is constructed from a Gauss diagram of a knotoid $K$.

\begin{definition} \cite{FLV} \label{def3.6} {\rm 
	A \textit{Gauss diagram} $G(D)$ of a diagram $D$ of a knotoid $K$ consists of an oriented arc of a circle with chords connecting the preimages of the crossings of $D$ according to the following rules:
	\begin{itemize}
		\item each chord corresponds to a crossing in the knotoid diagram;
		\item each chord is oriented from the overpassing branch to the underpassing branch;
		\item the sign of the crossing $\operatorname{sgn} (c)$ is indicated near the tail of the corresponding chord, see Fig.~\ref{fig8}.
	\end{itemize}
}
\begin{figure}[htbp] 
	\centering
	\begin{tikzpicture}[scale=0.6]
		\draw[black, thick] (0,0) arc (30:330:2);
		\filldraw[black] (0,0) circle (.1);
		\filldraw[black] (0,-2) circle (.1);
		\draw [black, very thick,->](-1.7,1) -- (-1.7,-3);
		\node [above] at(-1.5,1){$\operatorname{sgn}(c)$};
		\node[right] at (-1.5,-1){$c$};
	\end{tikzpicture}
	\caption{A chord of a Gauss diagram corresponding to a crossing $c$.}	\label{fig8}
\end{figure}
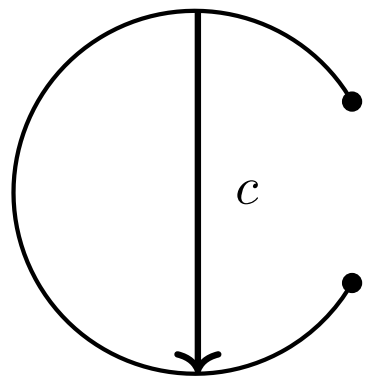
\end{definition}

To construct $H_D(t,y,z)$, we need the following notions and notation:
\begin{enumerate}
\item denote by $C(G(D))$ the set of chords of the Gauss diagram $G(D)$;
\item denote by $r^{+}(c)$ (resp., $r^-(c)$) the set of chords in $G(D) \setminus \{c\}$ that intersect the chord $c$ from left to right with respect to the orientation of $c$ and have sign $+$ (resp., $-$);
\item denote by $\ell^+(c)$ (resp., $\ell^-(c)$) the set of chords in $G(D) \setminus \{c\}$ that intersect the chord $c$ from right to left with respect to the orientation of $c$ and have sign $+$ (resp., $-$);
\item the \textit{degree} $d(c)$ of a chord $c$ is defined as the integer
$$
d(c) = |r^+(c)|-|r^-(c)| - |\ell^+(c)|+|\ell^-(c)|; 
$$
\item for $n \in \mathbb{N}$, set
$$
C_n(c) = \{e \in C(G(D))|\ e\cap c \neq \emptyset,\  \gcd (d(c), d(e)) = n \},
$$
where by definition we set $\gcd(m, 0) = |m|$ and $\gcd(0, 0) = 0$; 
\item set $r(c) = r^+(c)\cup r^-(c)$ and $\ell(c) = \ell^+(c)\cup \ell^-(c)$; 
\item set $r^n(c) = r(c) \cap C_n(c)$ and $\ell^n(c) = \ell(c)\cap C_n(c)$; 
\item define $\Phi_c(k) = k \mod |d(c)|$ for $k \in \mathbb{Z}$; by definition, set $\Phi_c(k) = k$ when $d(c) = 0$.
\end{enumerate}

\begin{definition} \label{def3.7} {\rm 
	Let $c \in C(G(D))$ and $n \in \mathbb{N}$. The following function is called the \textit{$n$-index function}:}
\begin{equation} \label{eqn:Ind}
	\operatorname{Ind}_c^n(z) := \displaystyle\sum_{e \in r^n(c)} \operatorname{sgn}(e) \cdot z^{\Phi_c(d(e))} - \displaystyle\sum_{e \in \ell^n(c)}\operatorname{sgn}(e)\cdot z^{\Phi_c(-d(e))}.
\end{equation}
\end{definition}

Note that $\operatorname{Ind}_c^n(z)$ takes values in $\mathbb Z[z^{\pm1}] / (z^{|d(c)|} - 1)$. 

We now define the three-variable function introduced in~\cite{FLV}.

\begin{definition} \cite{FLV} \label{3.8} {\rm 
	Let $D$ be a knotoid diagram. Set}
\begin{equation} \label{eqn:H}
	H_D(t,y,z) := \displaystyle \sum_{c \in C(G(D)), \ n \in \mathbb{N}} \operatorname{sgn}(c) \cdot (t^{\operatorname{Ind}_c^n(z)} - 1) \cdot y^n.
\end{equation}
\end{definition}

Observe that in the connected sum of two knotoids, in the Gauss diagram we first traverse all crossings of the first knotoid, and then all crossings of the second knotoid. Consequently, the chords of the first knotoid do not intersect the chords of the second knotoid. Thus, $H_{K\#K'}(t,y,z)=H_K(t,y,z)+H_{K'}(t,y,z)$.

\begin{definition} \label{def3.9} {\rm 
	The \textit{$y$-span} $\mathbf{span}_y(H_K(t,y,z))$ is defined as the difference between the maximum and minimum degrees of $y$ in $H_K(t,y,z)$.
	If $H_K(t,y,z)$ is the identically zero function then we set $\mathbf{span}_y(H_K(t,y,z))=0$.
	}
	\end{definition}
	
	\begin{definition} \label{def3.11} {\rm 
	The \textit{$z$-span} $\mathbf{span}_z(H_K(t,y,z))$ is defined as the difference between the maximum and minimum degrees of $z$ in $H_K(t,y,z)$.
		If $H_K(t,y,z)$ is the identically zero function then we set $\mathbf{span}_z(H_K(t,y,z))=0$.
	
	}
	\end{definition}
	
	Note that $\mathbf{span}_y(H_{K}(t,y,z))$ and $\mathbf{span}_z(H_{K}(t,y,z))$ are invariants of a spherical knotoid $K$, since $H_K(t,y,z)$ is an invariant of $K$.
	
	In~\cite{FLV}, the following property was established.
	
\begin{theorem} \cite[Th~1.4]{FLV} \label{theorem3.4}
	Let $K$ and $K'$ be two planar knotoids such that $K$ can be transformed into $K'$ by a finite sequence of crossing changes. Then
	$$
	H_K(t,y,z) - H_{K'}(t,y,z) = \displaystyle \sum_{n \in \mathbb{N}} \left( \sum_{m \in \mathbb{N}} a_{n_m} \left( t^{z_n^m} + t^{-(z^{-1})_n^m} - 2 \right) \right) y^n,
	$$
	\text{where}\ $a_{n_m} \in \mathbb{Z}$ and $z_n^m$ is some polynomial of the variable z. 
\end{theorem}

We note a property of this invariant for spherical knotoids.  

\begin{corollary}\label{lemma3.5}
	If $K$ is a spherical knotoid, then $H_K(t,y,z) = H_K(t^{-1},y,z^{-1})$.
\end{corollary}

\begin{proof}
	In~\cite[Pr.~3.6]{GuKa}, it was shown that any spherical knotoid can be transformed into the trivial knotoid $U$ by a finite sequence of crossing changes. Since $H_U(t,y,z) = 0$, by Theorem~\ref{theorem3.4} we have
	$$
	H_K(t,y,z) = \displaystyle \sum_{n \in \mathbb{N}} \left( \sum_{m \in \mathbb{N}} a_{n_m} \left( t^{z_n^m} + t^{-(z^{-1})_n^m} - 2 \right) \right) y^n, 
	$$ 
	whence
	$$
	H_K(t^{-1},y,z^{-1}) = \displaystyle \sum_{n \in \mathbb{N}} \left( \sum_{m \in \mathbb{N}} a_{n_m} \left(t^{-(z^{-1})_n^m} + t^{z_n^m} - 2 \right) \right) y^n, 
	$$ 
	and thus $H_K(t,y,z) = H_K(t^{-1},y,z^{-1})$. 
\end{proof}

\begin{corollary}\label{lemma3.7}
	For any spherical knotoid $K$, the number $\mathbf{span}_z(H_K(t, y, z))$ is even and equals twice the maximum degree of the variable $z$.
\end{corollary}

\section{The family of spiral knotoids $S_m$} \label{sec5}

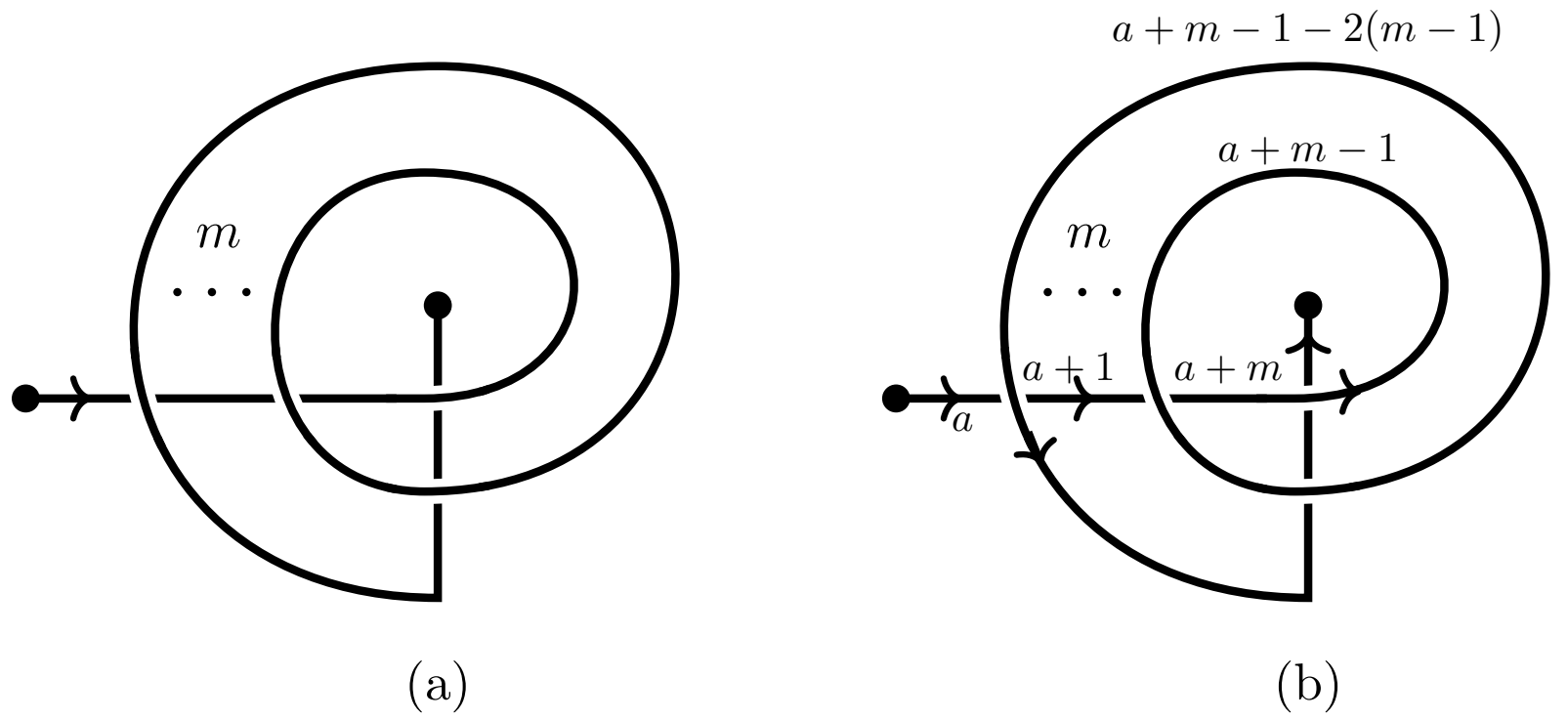
\begin{figure}[htbp]
	\begin{center}
		\begin{tikzpicture}[scale=0.9]
			\begin{knot}[consider self intersections, ignore endpoint intersections=false, clip radius=8pt]
				\strand[ultra thick,black]
				(0,1.5) -- (3,1.5)
				.. controls (4.5,1.5) and (4.5,3.2) .. (3,3.2)
				.. controls (1.5,3.2) and (1.5,0.8) .. (3,0.8)
				.. controls (5.5,0.8) and (5.5,4.) .. (3.1,4.)
				.. controls (0.,4.) and (0.1,0) .. (3.1,0)
				-- (3.1,2.2);
				\flipcrossings{1, 2}
			\end{knot}
			\draw [black, very thick,->](0,1.5) -- (0.5,1.5);
			\node[black, scale=1.5] at (1.45, 2.3) {$\dots$};
			\node[black, above] at (1.45, 2.5) {$m$};
			\fill[black] (0,1.5) circle (3pt);
			\fill[black] (3.1,2.2) circle (3pt);
			\node[black, above] at (3.1,-1.0) {(a)};
		\end{tikzpicture}
		\qquad
		\begin{tikzpicture}[scale=0.9]
			\begin{knot}[consider self intersections, ignore endpoint intersections=false, clip radius=8pt]
				\strand[ultra thick,black]
				(0,1.5) -- (3,1.5)
				.. controls (4.5,1.5) and (4.5,3.2) .. (3,3.2)
				.. controls (1.5,3.2) and (1.5,0.8) .. (3,0.8)
				.. controls (5.5,0.8) and (5.5,4.) .. (3.1,4.)
				.. controls (0.,4.) and (0.1,0) .. (3.1,0)
				-- (3.1,2.2);
				\flipcrossings{1, 2}
			\end{knot}
			\draw [black, very thick,->](0,1.5) -- (0.5,1.5);
			\node[black, scale=1.5] at (1.45, 2.3) {$\dots$};
			\node[black, above] at (1.45, 2.5) {$m$};
			\fill[black] (0,1.5) circle (3pt);
			\fill[black] (3.1,2.2) circle (3pt);
			\draw [black, very thick,->](1,1.5) -- (1.5,1.5);
			\draw [black, very thick,->](3.45,1.55) -- (3.5,1.55);
			\draw [black, very thick,->](3.1,1.75) -- (3.1,2);
			\draw [black, very thick,->](1,1.25) -- (1.1,1);
			\node[below, scale=0.8] at (0.5, 1.5) {$a$};
			\node[above, scale=0.8] at (1.3, 1.5) {$a+1$};
			\node[above, scale=0.8] at (2.5, 1.5) {$a+m$};
			\node[below, scale=0.8] at (3.1, 3.6) {$a+m-1$};
			\node[above, scale=0.8] at (3.1, 4.) {$a+m-1-2(m-1)$};
			\node[black, above] at (3.1,-1.0) {(b)};
		\end{tikzpicture}
	\end{center}
	\caption{The spiral knotoid $S_m$ and labels on its arcs.} \label{fig6}
\end{figure}

Denote by $S_m$, $m \geq 1$, the spiral knotoid in the $2$-sphere shown in Fig.~\ref{fig6}(a), where the number of crossings is equal to $2m$.

\begin{figure}[htbp]
	\begin{center}
		\begin{tikzpicture}[scale=0.9]
			\begin{knot}[consider self intersections, ignore endpoint intersections=false, clip radius=8pt]
				\strand[ultra thick,black]
				(0,1.5) -- (3,1.5)
				.. controls (4.5,1.5) and (4.5,3.2) .. (3,3.2)
				.. controls (1.5,3.2) and (1.5,0.8) .. (3,0.8)
				.. controls (5.5,0.8) and (5.5,4.) .. (3.1,4.)
				.. controls (0.,4.) and (0.1,0) .. (3.1,0)
				-- (3.1,2.5);
				\flipcrossings{1, 2}
			\end{knot}
			\draw [black, very thick,->](0,1.5) -- (0.5,1.5);
			\node[black, scale=1.5] at (1.45, 2.3) {$\dots$};
			\node[black, above] at (1.45, 2.5) {$m$};
			\fill[black] (0,1.5) circle (3pt);
			\fill[black] (3.1,2.5) circle (3pt);
			\draw[ultra thick, red, dashed]  (3.55,2) circle[x radius=0.7, y radius=0.3];
			\draw[->] (5., 1.5) -- (6., 1.5) node[midway, above] {$H(2)$};
			\begin{knot}[consider self intersections, ignore endpoint intersections=false, clip radius=8pt]
				\strand[ultra thick,black]
				(7,1.5) -- (10,1.5)
				.. controls (11.5,1.5) and (11.5,3.2) .. (10,3.2)
				.. controls (8.5,3.2) and (8.5,0.8) .. (10,0.8)
				.. controls (12.5,0.8) and (12.5,4.) .. (10.1,4.)
				.. controls (7.1,4.) and (7.1,0) .. (10.1,0)
				-- (10.1,2.5);
				\flipcrossings{1, 2}
			\end{knot}
			\draw [black, very thick,->](7,1.5) -- (7.5,1.5);
			\node[black, scale=1.5] at (8.45, 2.3) {$\dots$};
			\node[black, above] at (8.45, 2.5) {$m$};
			\fill[black] (7,1.5) circle (3pt);
			\fill[black] (10.1,2.5) circle (3pt);
			\fill[white] (10.55,2) circle[x radius=0.7, y radius=0.3];
			\draw[ultra thick, black]  (10.1,2.2) arc (-180:0: 0.5 and 0.15);
			\draw[ultra thick, black]  (10.1,1.7) arc (180:0: 0.34 and 0.15);
			\draw[ultra thick, red, dashed]  (10.55,2) circle[x radius=0.7, y radius=0.3];
		\end{tikzpicture}
	\end{center}
	\caption{An $H(2)$-move unknotting the knotoid $S_m$.} \label{fig7}
\end{figure}
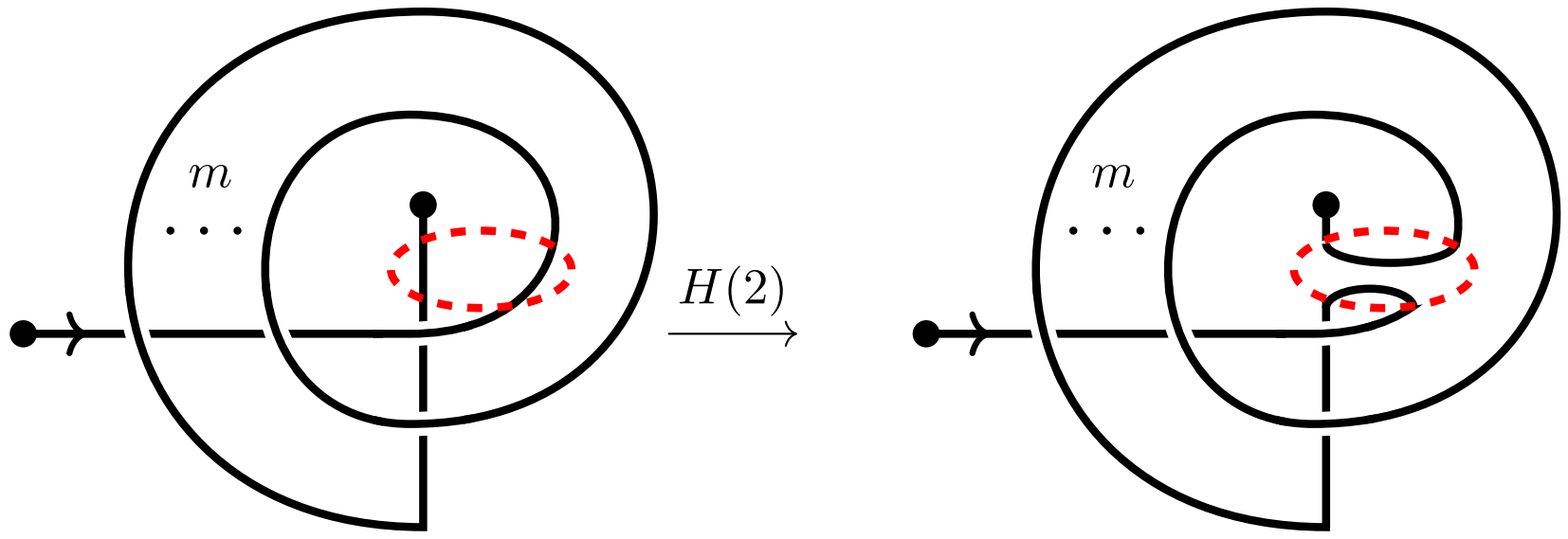

\begin{lemma} \label{lemma3.11}
	For any $n \geq 2$ and any $m \geq 1$, the spiral knotoid $S_m$ can be transformed into the trivial knotoid $U$ by a single $H(n)$-move.
\end{lemma}

\begin{proof}
	It is easy to see that the $H(2)$-move shown in Fig.~\ref{fig7} unknots $S_m$. In~\cite[Th.~5]{JLM}, the inequality
	$$
	d_{H(n)}(K,K')\geq d_{H(n+1)}(K,K')
	$$
	was proved for knots. Since the proof uses only local $H(n)$-moves, it carries over verbatim to knotoids. Hence the inequality also holds for knotoids. Thus, $S_m$ can be transformed into the trivial knotoid $U$ by a single $H(n)$-move for any $n\geq 2$. 
\end{proof}

\begin{lemma} \label{lemma3.12}
	For any $p, q \geq 1$, the connected sum of spiral knotoids $S_p \# S_q$ can be transformed into the knotoid $S_p$ by a single $H(n)$-move.
\end{lemma}

\begin{proof}
	By Lemma~\ref{lemma3.11}, the summand $S_q$ in the connected sum $S_p \# S_q$ can be transformed into the trivial knotoid $U$ by a single $H(n)$-move without affecting the summand $S_p$.
\end{proof}

We now compute the invariants of the spiral knotoids $S_m$. 

\begin{lemma} \label{lemma3.13}
	For the spherical spiral knotoid $S_m$, $m \geq 1$, we have $h(S_m) = m$.
\end{lemma}

\begin{proof}
	Assign integer labels to the arcs of the knotoid diagram following its orientation, starting with the label $a$ on the arc incident to the tail of the knotoid, see Fig.~\ref{fig6}(b). As seen from the diagram, the sign of every crossing is $+1$. Following the orientation of the knotoid, we pass the first $m$ crossings from right to left (hence the label increases by $1$ at each step). Then we pass the $(m+1)$-st crossing from left to right. Next, passing each of the $m-1$ turns of the spiral changes the label by $-2$. We obtain the labels shown in Fig.~\ref{fig6}(b). Let $c_1$ be the first crossing in the diagram of $S_m$. Then $w(c_1) = \operatorname{sgn}(c_1)\cdot \left[ (a+m-1-2(m-1))-(a+1) \right] = (+1) \cdot \left[ (a+1-m)-(a+1) \right] = -m$. By Theorem~\ref{theorem3.3}, there exists a crossing $c$ such that $w(c) = m$. Since all crossing signs of the knotoid are positive, the monomial $t^{m}$ does not cancel when computing $P_{S_m}(t)$. By Theorem~\ref{theorem3.2}, we obtain $h(S_m) \geq m$, while according to the diagram in Fig.~\ref{fig6}(a), $h(S_m) \leq m$. Consequently, $h(S_m) = m$.
\end{proof}	

\begin{corollary} \label{lemma3.8P}
	(1) For the spherical spiral knotoid $S_m$, $m \geq 1$, we have $\mathbf{span}(P_{S_m}(t)) = 2m$. \\ 
	(2) For the connected sum of spherical spiral knotoids, we have  $\mathbf{span}(P_{S_p\#S_q}) = \max\{2p,2q\}$. 
\end{corollary}

\begin{proof}
	(1) As shown in the proof of Lemma~\ref{lemma3.13}, the maximum degree of $P_{S_m}(t)$ is $m$. By Theorem~\ref{theorem3.3}, we obtain $\mathbf{span}(P_{S_m}(t)) = 2m$.  \\
	(2) This follows from the equality $P_{S_{p}\#S_{q}} = P_{S_{p}} + P_{S_{q}}$. 
\end{proof}

Now we compute the function $H_K(t, y, z)$ for $S_m$. 

\begin{lemma} \label{lemma3.3}
	For the spherical spiral knotoid $S_m$, $m \geq 2$, we have $\mathbf{span}_y(H_{S_m}(t,y,z)) = m-1$.
\end{lemma}

\begin{proof}
	Label the crossings of the diagram of $S_m$ as shown in Fig.~\ref{fig9}(a).
	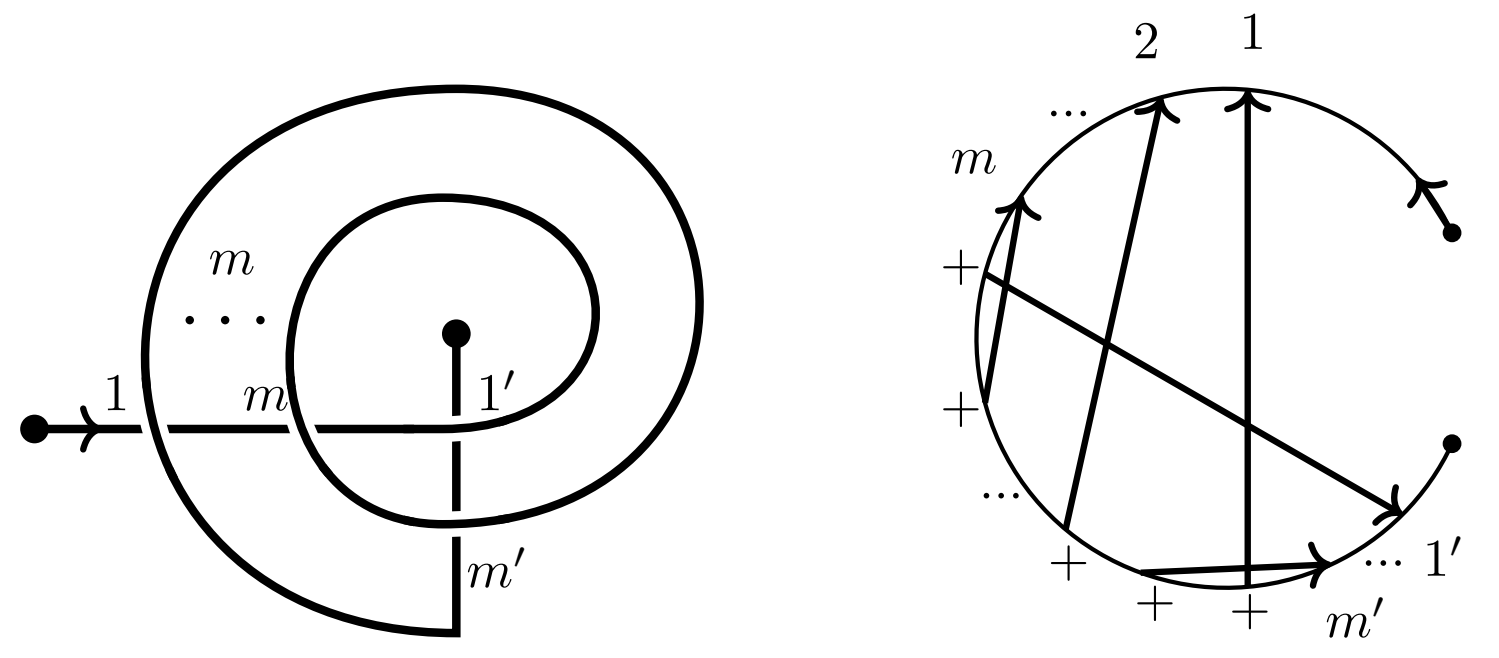
\begin{figure}[htbp]
		\begin{center}
			\begin{tikzpicture}[scale=0.9]
				\begin{knot}[consider self intersections, ignore endpoint intersections=false, clip radius=8pt]
					\strand[ultra thick,black]
					(0,1.5) -- (3,1.5)
					.. controls (4.5,1.5) and (4.5,3.2) .. (3,3.2)
					.. controls (1.5,3.2) and (1.5,0.8) .. (3,0.8)
					.. controls (5.5,0.8) and (5.5,4.) .. (3.1,4.)
					.. controls (0.,4.) and (0.1,0) .. (3.1,0)
					-- (3.1,2.2);
					\flipcrossings{1, 2}
				\end{knot}
				\draw [black, very thick,->](0,1.5) -- (0.5,1.5);
				\node[black, scale=1.5] at (1.45, 2.3) {$\dots$};
				\node[black, above] at (1.45, 2.5) {$m$};
				\node[above] at (0.6, 1.5) {$1$};
				\node[above] at (1.7, 1.5) {$m$};
				\node[above] at (3.4, 1.5) {$1'$};
				\node[above] at (3.4, 0.2) {$m'$};
				\fill[black] (0,1.5) circle (3pt);
				\fill[black] (3.1,2.2) circle (3pt);
				\node[above] at (3.1, -1.5) {(a)};
			\end{tikzpicture}
			\qquad
			\begin{tikzpicture}[scale=0.55]
				\draw[black, thick] (25:3) arc (25:335:3);
				\filldraw[black] (25:3) circle (.1); node[right=5pt] {l};
				\filldraw[black] (335:3) circle (.1); node[right=5pt] {h};
				\draw [black, very thick,->] (25:3) -- (40:3);
				\fill (85:3) coordinate (c'4) circle (0.5pt);
				\node at (85:3.3) {$ $};
				\node at (85:3.7) {$1$};
				\fill (105:3) coordinate (c3) circle (0.5pt);
				\node at (105:3.3) {$ $};
				\node at (105:3.7) {$2$};
				\node at (125:3.3) {$...$};
				\fill (145:3) coordinate (c5) circle (0.5pt);
				\node at (145:3.3) {$ $};
				\node at (145:3.7) {$m$};
				\fill (165:3) coordinate (c2) circle (0.5pt);
				\node at (165:3.3) {$+$};
				\node at (165:3.7) {$ $};
				\fill (195:3) coordinate (c'5) circle (0.5pt);
				\node at (195:3.3) {$+$};
				\node at (195:3.7) {$ $};
				\node at (215:3.3) {$...$};
				\fill (230:3) coordinate (c'3) circle (0.5pt);
				\node at (235:3.3) {$+$};
				\node at (230:3.7) {$ $};
				\fill (250:3) coordinate (c'6) circle (0.5pt);
				\node at (255:3.3) {$+$};
				\node at (250:3.7) {$$};
				\fill (275:3) coordinate (c4) circle (0.5pt);
				\node at (275:3.3) {$+$};
				\node at (275:3.7) {$ $};
				\fill (295:3) coordinate (c6) circle (0.5pt);
				\node at (295:3.3) {$$};
				\node at (295:3.7) {$m'$};
				\node at (305:3.3) {$...$};
				\fill (315:3) coordinate (c'2) circle (0.5pt);
				\node at (315:3.3) {};
				\node at (315:3.7) {$1'$};
				\draw[black, very thick, ->] (c2) -- (c'2);
				\draw[black, very thick, ->] (c'3) -- (c3);
				\draw[black, very thick, ->] (c4) -- (c'4);
				\draw[black, very thick, ->] (c'5) -- (c5);
				\draw[black, very thick, ->] (c'6) -- (c6);
				\node[above] at (0, -6) {(b)};
			\end{tikzpicture}
		\end{center}
		\caption{Labeling of crossings and the Gauss diagram for the knotoid $S_m$.} \label{fig9}
	\end{figure}
	Namely, moving from the tail of the knotoid to its head, we pass each crossing twice in the following order:
	$$
	\{ 1, 2, \ldots, (m-1), m, 1', m, 2', (m-1), \ldots, m', 1, m', (m-1)', \ldots, 1' \}.
	$$
	Note that the crossings $\{ 1, 2, \ldots, m\}$ are first passed as under-crossings and then as over-crossings. Conversely, the crossings $\{ 1', 2', \ldots, m' \}$ are first passed as over-crossings and then as under-crossings. Taking into account the orientation of the knotoid, we obtain $\operatorname{sgn}(k) = +1$ and $\operatorname{sgn}(k') = +1$ for $k = 1, 2, \ldots, m$.
	
	From the Gauss diagram of $S_m$ shown in Fig.~\ref{fig9}(b), we see that the chords in the set $M = \{ 1, 2, \ldots, m\}$ are pairwise disjoint. Similarly, the chords in the set $M' = \{ 1', 2', \ldots, m' \}$ are pairwise disjoint. If $k' \in M'$ intersects $s \in M$, then $k'$ is oriented from left to right with respect to $s$. Similarly, if $k\in M$ intersects $s' \in M'$, then $k$ is oriented from right to left with respect to $s'$. Thus, we obtain the sets
	\begin{equation}
		\begin{gathered}
			r^+(1) = \{1', \ldots, m'\}, \quad  r^+(2) = \{1', \ldots, (m-1)'\},  \quad \ldots, \quad r^+(m) = \{1'\}; \cr
			\ell^+(1') = \{1, \ldots, m\},  \quad  l^+(2') = \{1, \ldots, (m-1)\}, \quad \ldots, \quad l^+(m') = \{1\}
		\end{gathered} \label{eqn:sets}
	\end{equation}
	and the values of the chord degrees:
	\begin{equation}
		d(1) = m, \quad \ldots, \quad d(m) = 1, \quad d(1') = -m, \quad \ldots, \quad d(m') = -1. \label{eqn:degrees}
	\end{equation}
	By definition, the degree $m$ of the variable $y$ can appear at a chord $e_1$ only if there exists a chord $e_2$ such that $\gcd (d(e_1),d(e_2))=m$. In the case of the knotoid $S_m$, only the chords $1\in M$ and $1' \in M'$ have this property. Thus, the maximum degree of the variable $y$ in $H_{S_m} (t, y, z)$ is equal to $m$.
	
	Since the degree of the chord $m \in M$ is $d(m) = 1$, its greatest common divisor with the degree of any other chord is $1$. Consequently, the minimum degree of the variable $y$ in $H_{S_m} (t, y, z)$ is equal to~$1$. Since the signs of all chords are positive, collecting like terms in~(\ref{eqn:H}) when computing $H_{S_m}(t, y, z)$ does not affect the difference between the maximum and minimum degrees of the variable $y$. Thus, $\mathbf{span}_y(H_{S_m}(t,y,z)) = m-1$.
\end{proof}

\begin{lemma} \label{lemma3.6}
	For the spherical spiral knotoid $S_m$, $m \geq 2$, we have $\mathbf{span}_z(H_{S_m}(t,y,z)) = 2(m-1)$.
\end{lemma}

\begin{proof}
	As in Lemma~\ref{lemma3.3}, following Fig.~\ref{fig9}, we obtain the sets~(\ref{eqn:sets}) and the values of the chord degrees~(\ref{eqn:degrees}). When computing the function $\operatorname{Ind}_c^m(z)$ by formula~(\ref{eqn:Ind}), we see that the maximum possible degree of the variable $z$ is $m-1$, since the largest value of $d(c)$ is $m$ and hence the largest remainder upon division (the value of the function $\Phi$) is $m-1$. Let us show that such a degree actually occurs in one of the terms.
	
	Consider the chord $1 \in M$. The set $r^+(1)$ contains the chord $2' \in M'$, with $d(1) = m$ and $d(2') = -(m-1)$. Then $2' \in r^1(1)$ and $\Phi_1(d(2')) = m-1$. Since the signs of all chords are positive, collecting like terms when computing $H_{S_m}(t, y, z)$ does not affect the difference between the maximum and minimum degrees of the variable $z$. Thus, $\mathbf{span}_z(H_{S_{m}}(t,y,z)) = 2(m-1)$.
\end{proof}

\section{Quotients of the Gordian graph of knotoids by invariants} \label{sec6}

\begin{theorem}\label{theorem4.1}
	Let $G(H(n), h)$, $G(H(n), \mathbf{span}(P(t)))$, $G(H(n), \mathbf{span}_y(H_K(t, y, z)))$, and $G(H(n),$ $\mathbf{span}_z(H_K(t, y, z)))$ be the quotient graphs of the $H(n)$-Gordian graph of knotoids $G(H(n))$ by the height $h$, the span $\mathbf{span}(P(t))$ of the affine index polynomial, the $y$-span $\mathbf{span}_y(H_K(t, y, z))$ and the $z$-span $\mathbf{span}_z(H_K(t, y, z))$ of the invariant $H_K(t, y, z)$, respectively. Then each of these quotient graphs is isomorphic to a complete graph on countably infinitely many vertices and is $0$-hyperbolic in the sense of Gromov.  
\end{theorem}

\begin{proof}
	We consider the four quotient graphs sequentially. 
	
	1. The graph $G(H(n), h)$. 
	By Theorem~\ref{theorem3.1} and Lemma~\ref{lemma3.13}, for spiral knotoids we obtain $h(S_p \# S_q) = h(S_p) + h(S_q) = p+q$, while by Lemma~\ref{lemma3.12}, $d_{H(n)}(S_p, S_p \# S_q) = 1$, since the summand $S_q$ can be transformed into the trivial knotoid $U$ by a single $H(n)$-move. Therefore, for any two classes $[K]$ and $[K']$ in the quotient of $G(H(n))$ by the height $h$, there exist representatives connected by a single $H(n)$-move. Thus, the graph $G(H(n), h)$ is complete. Inequality~(\ref{eqn0}) for the distances between any four vertices holds with $\delta = 0$, meaning that the graph is $0$-hyperbolic.
	
	2. The graph $G(H(n), \mathbf{span}(P(t)))$.
	By Corollary~\ref{lemma3.8P}, $\mathbf{span}(P_{S_{p}\# S_{q}}) = \max\{2p, 2q\}$. By Lemma~\ref{lemma3.12}, $d_{H(n)}(S_p, S_p \# S_q) = 1$. Therefore, for any two classes $[K]$ and $[K']$ in the quotient by $\mathbf{span}(P)$, there exist representatives connected by a single $H(n)$-move. Thus, the graph $G(H(n), \mathbf{span}(P))$ is complete. Inequality~(\ref{eqn0}) for the distances between any four vertices holds with $\delta = 0$, meaning that the graph is $0$-hyperbolic.
	
	3. The graph $G(H(n), \mathbf{span}_y(H_K(t, y, z)))$.
	By Lemma~\ref{lemma3.3}, the set of values of $\mathbf{span}_y(H_K(t, y, z))$ coincides with the set $\mathbb N \cup \{0\}$. From the additivity $H_{K \# K'}(t, y, z) = H_K(t, y, z) + H_{K'}(t, y, z)$, the positivity of the crossing signs of the spiral knotoids $S_m$, and Lemma~\ref{lemma3.3}, it follows that $\mathbf{span}_y(H_{S_{p+1} \# S_{q+1}}(t, y, z)) = \max\{p, q\}$. By Lemma~\ref{lemma3.12}, $d_{H(n)}(S_{p+1}, S_{p+1} \# S_{q+1}) = 1$. Therefore, for any two classes $[K]$ and $[K']$ in the quotient by $\mathbf{span}_y(H_K(t, y, z))$, there exist representatives connected by a single $H(n)$-move. Thus, the graph $G(H(n), \mathbf{span}_y(H_K(t, y, z)))$ is complete. Inequality~(\ref{eqn0}) for the distances between vertices holds with $\delta = 0$, meaning that the graph is $0$-hyperbolic.
	
	4. The graph $G(H(n), \mathbf{span}_z(H_K(t, y, z)))$. 
	By Lemma~\ref{lemma3.6}, the set of values of $\mathbf{span}_z(H_K(t, y, z))$ coincides with the set of all even positive integers and zero. From the additivity $H_{K \# K'} (t, y, z)= H_K(t, y, z) + H_{K'}(t, y, z)$, the positivity of the crossing signs of the spiral knotoids $S_m$, and Lemma~\ref{lemma3.6}, it follows that $\mathbf{span}_z(H_{S_{p+1}\# S_{q+1}}(t, y, z)) = \max\{2p,2q\} = 2p$ for $p > q$. By Lemma~\ref{lemma3.12}, \linebreak $d_{H(n)}(S_{p+1}, S_{p+1} \# S_{q+1}) = 1$. Therefore, for any two classes $[K]$ and $[K']$ in the quotient by $\mathbf{span}_z(H_K(t, y, z))$, there exist representatives connected by a single $H(n)$-move. Thus, the graph $G(H(n), \mathbf{span}_z(H_K(t, y, z)))$ is complete. Inequality~(\ref{eqn0}) for the distances between vertices holds with $\delta = 0$, meaning that the graph is $0$-hyperbolic.
	The theorem is proved.
\end{proof}

\section{Quotient of the Gordian graph of knotoids by the span of the Jones polynomial}   \label{sec7}

\begin{definition} \label{def3.13} {\rm 
		An \textit{oriented multi-knotoid diagram} in $S^2$ is an immersion into $S^2$ of an oriented interval $[0,1]$ and a finite number of oriented circles, such that all double points are endowed with over/under-crossing information, and the images of the endpoints of the interval are distinct from each other and from the double points.}
\end{definition}

\begin{definition} \label{def3.14} {\rm 
	An \textit{oriented spherical multi-knotoid} is an equivalence class of oriented multi-knotoid diagrams in $S^2$ up to Reidemeister moves applied in regions disjoint from the endpoints. }
\end{definition}

In~\cite{Tu}, by analogy with the bracket polynomial of links, the bracket polynomial for spherical multi-knotoids and then the normalized bracket polynomial were defined. 

\begin{definition} \cite{Tu} \label{def3.9bis} {\rm 
	For a diagram $D$ of a multi-knotoid in the sphere, the \textit{bracket polynomial} is defined as
	$$
	\langle  D \rangle = \displaystyle \sum_{s \in S} A^{\sigma_s}\left( -A^2-A^{-2}\right)^{|s|-1},
	$$
	where $S$ is the set of all states of the multi-knotoid diagram, $\sigma_s = n_+ - n_-$, where $n_+$ is the number of $A$-smoothings and $n_-$ is the number of $A^{-1}$-smoothings (see~(\ref{eqn7})) used to obtain the state $s$, and $|s|$ is the number of connected components of the state $s$. }
\end{definition}

\begin{definition} \cite{Tu} \label{def3.8bis} {\rm
	For a diagram $D$ of a multi-knotoid in the sphere, the \textit{normalized bracket polynomial} is defined as
	$$
	\langle  D \rangle_\circ = \displaystyle (-A^3)^{-wr(D)} \langle  D \rangle,
	$$
	where $wr(D) = \displaystyle \sum_{c \in C} \operatorname{sgn} (c)$ is the writhe and $C$ is the set of all crossings of the diagram $D$.}
\end{definition}

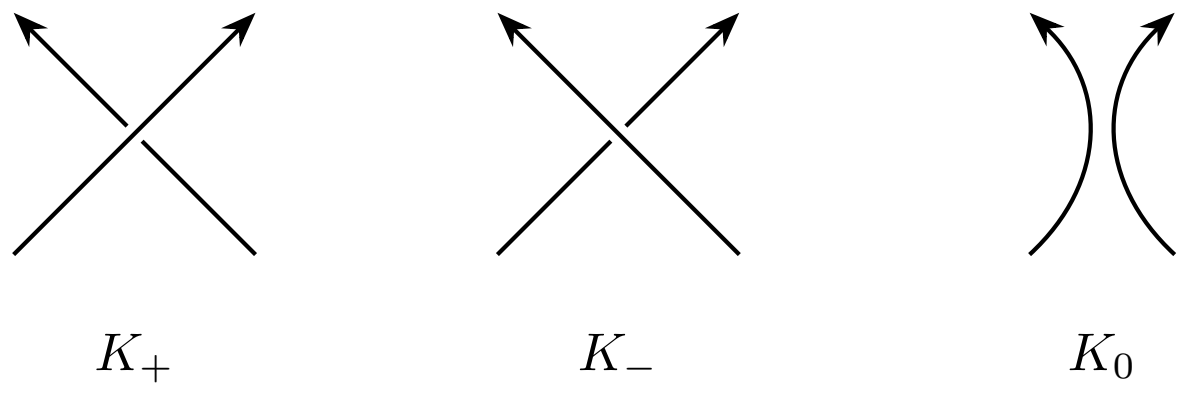
\begin{figure}[ht]
	\begin{tikzpicture}[baseline=-0.6ex, scale=0.8, >=Stealth]
		\begin{scope}
			\draw[thick, ->] (1,-1) -- (-1,1);
			\draw[white, line width=4pt] (-1,-1) -- (1,1);
			\draw[thick, ->] (-1,-1) -- (1,1);
			\node[below] at (0,-1.5) {$K_+$};
		\end{scope}
		\begin{scope}[shift={(4,0)}]
			\draw[thick, ->] (-1,-1) -- (1,1);
			\draw[white, line width=4pt] (1,-1) -- (-1,1);
			\draw[thick, ->] (1,-1) -- (-1,1);
			\node[below] at (0,-1.5) {$K_-$};
		\end{scope}
		\begin{scope}[shift={(8,0)}]
			\draw[thick, ->] (-0.6,-1)
			.. controls (0.05,-0.4) and (0.05,0.4) .. (-0.6,1);
			\draw[thick, ->] (0.6,-1)
			.. controls (-0.05,-0.4) and (-0.05,0.4) .. (0.6,1);
			\node[below] at (0,-1.5) {$K_0$};
		\end{scope}
	\end{tikzpicture}
	\caption{A Conway triple.} \label{fig14}
\end{figure}

In~\cite[Sec.~7.1]{Tu}, it was shown that $\langle  D \rangle_\circ = \displaystyle (-A^3)^{-wr(D)} \langle  D \rangle$ is invariant under Reidemeister moves. That is, the normalized bracket polynomial is an invariant of spherical multi-knotoids. Moreover (\cite[Sec.~8.4]{Tu}), the normalized bracket polynomial satisfies the skein relation
\begin{equation}\label{eqn6}
	-A^4\langle K_+ \rangle_\circ + A^{-4}\langle K_- \rangle_\circ = (A^2-A^{-2})\langle K_0 \rangle_\circ,
\end{equation}
where $(K_+,K_-,K_0)$ is a Conway triple, see Fig.~\ref{fig14}.

Recall that the bracket polynomials of an initial diagram and the diagrams obtained by performing an $A$-smoothing and an $A^{-1}$-smoothing at one of its crossings are related as follows:
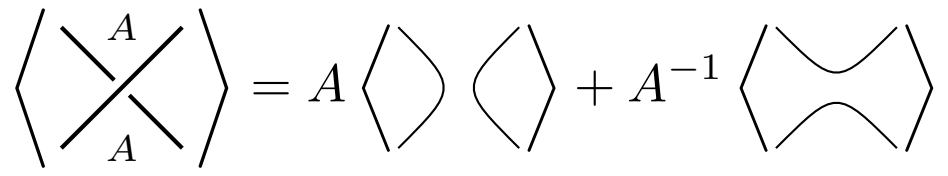
\begin{equation}\label{eqn7}
	\left\langle
	\begin{tikzpicture}[baseline=-0.6ex, scale=0.4]
		\draw[thick] (1,-1) -- (-1,1);
		\draw[white, line width=4pt] (-1,-1) -- (1,1);
		\draw[thick] (-1,-1) -- (1,1);
		\node[black] at (0, 1) {\scriptsize $A$};
		\node[black] at (0, -1) {\scriptsize $A$};
	\end{tikzpicture}
	\right\rangle = A \left\langle
	\begin{tikzpicture}[baseline=-0.6ex, scale=0.4]
		\draw[knot] (-1,-1) .. controls (0,0) .. (-1,1);
		\draw[knot] (1,-1) .. controls (0,0) .. (1,1);
	\end{tikzpicture}
	\right\rangle + A^{-1} \left\langle
	\begin{tikzpicture}[baseline=-0.6ex, scale=0.4]
		\draw[knot] (-1,-1) .. controls (0,0) .. (1,-1);
		\draw[knot] (-1,1) .. controls (0,0) .. (1,1);
	\end{tikzpicture}
	\right\rangle
\end{equation}
From~(\ref{eqn7}), we obtain the relation for normalized bracket polynomials: 
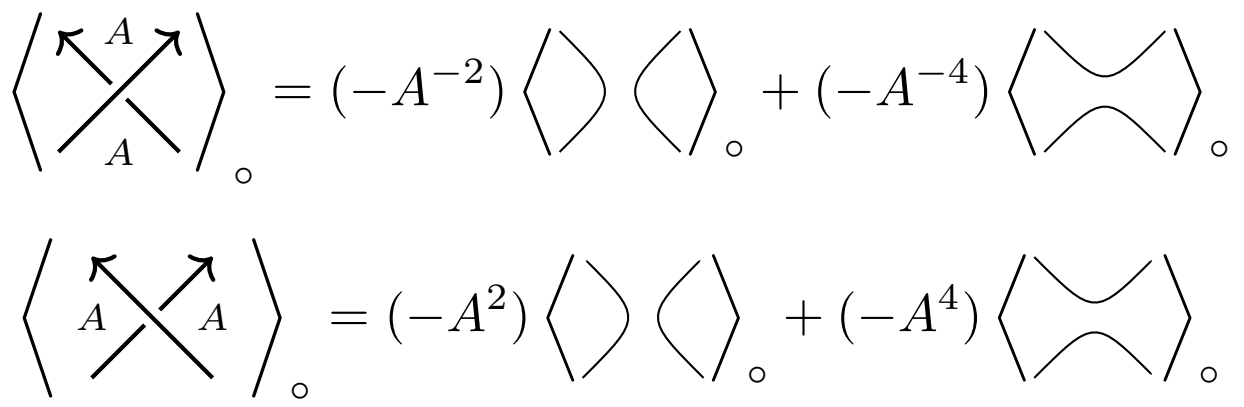
\begin{equation}\label{eqn8}
	\begin{gathered}
		\left\langle
		\begin{tikzpicture}[baseline=-0.6ex, scale=0.4]
			\draw[thick, ->] (1,-1) -- (-1,1);
			\draw[white, line width=4pt] (-1,-1) -- (1,1);
			\draw[thick, ->] (-1,-1) -- (1,1);
			\node[black] at (0, 1) {\scriptsize $A$};
			\node[black] at (0, -1) {\scriptsize $A$};
		\end{tikzpicture}
		\right\rangle_\circ = (-A^{-2}) \left\langle
		\begin{tikzpicture}[baseline=-0.6ex, scale=0.4]
			\draw[knot] (-1,-1) .. controls (0,0) .. (-1,1);
			\draw[knot] (1,-1) .. controls (0,0) .. (1,1);
		\end{tikzpicture}
		\right\rangle_\circ + (-A^{-4}) \left\langle
		\begin{tikzpicture}[baseline=-0.6ex, scale=0.4]
			\draw[knot] (-1,-1) .. controls (0,0) .. (1,-1);
			\draw[knot] (-1,1) .. controls (0,0) .. (1,1);
		\end{tikzpicture}
		\right\rangle_\circ  \\[6pt]
		\left\langle
		\begin{tikzpicture}[baseline=-0.6ex, scale=0.4]
			\draw[thick, ->] (-1,-1) -- (1,1);
			\draw[white, line width=4pt] (1,-1) -- (-1,1);
			\draw[thick, ->] (1,-1) -- (-1,1);
			\node[black] at (-1, 0) {\scriptsize $A$};
			\node[black] at (1, 0) {\scriptsize $A$};
		\end{tikzpicture}
		\right\rangle_\circ = (-A^2) \left\langle
		\begin{tikzpicture}[baseline=-0.6ex, scale=0.4]
			\draw[knot, ->] (-1,-1) .. controls (0,0) .. (-1,1);
			\draw[knot, ->] (1,-1) .. controls (0,0) .. (1,1);
		\end{tikzpicture}
		\right\rangle_\circ + (-A^4) \left\langle
		\begin{tikzpicture}[baseline=-0.6ex, scale=0.4]
			\draw[knot, ->] (-1,-1) .. controls (0,0) .. (1,-1);
			\draw[knot, ->] (-1,1) .. controls (0,0) .. (1,1);
		\end{tikzpicture}
		\right\rangle_\circ 
	\end{gathered}
\end{equation}

\begin{definition} \label{def3.12} {\rm 
	Let $K$ be a spherical multi-knotoid. The \textit{Jones polynomial} of $K$ is defined as $V_K(t) :=\langle K \rangle_\circ |_{A = t^{-1/4}}$. The \textit{span} $\mathbf{span}(V_K(t))$ of the Jones polynomial $V_K(t)$ is defined as the difference between the maximum and minimum degrees of $t$ in $V_K(t)$. }
\end{definition}

It is well known~\cite[Ch.~3]{Li} that for knots, the span of the Jones polynomial takes only integer values. However, this is not true for knotoids. For example, the span of the Jones polynomial for the knotoid $\phi$ from~\cite[Sec.~7.4]{Tu} is equal to ${{3}\over{2}}$. We show that for knotoids, the span of the Jones polynomial can take only integer and half-integer values.

\begin{proposition}\label{prop3.1}
	Let $K$ be a spherical knotoid. Then $\mathbf{span}(V_K(t)) \in \frac{1}{2} \mathbb Z$.
\end{proposition}

\begin{proof} 
	From the definition, it follows that $\mathbf{span}(V_K(t)) = \frac{1}{4}\mathbf{span}(\langle D \rangle_\circ (A))$. Let us show that $\mathbf{span}(\langle D \rangle_\circ (A))$ is divisible by $2$. Let $s$ be a state of a diagram of the knotoid $K$, and let $s'$ be a state differing from $s$ by the type of smoothing at a single crossing. Then, depending on the smoothing type, $\sigma_{s'} - \sigma_{s} = \pm 2$. If the number of connected components does not change under the smoothing, then the parity of the degree of $A$ does not change. If the number of connected components changes, then the degree of $A$ changes by an even integer due to the factor $\left( -A^2-A^{-2}\right)$. Thus, passing to a new state of a knotoid diagram preserves the parity of the degree of $A$. Hence, the maximum and minimum degrees of $A$ have the same parity, so the span of the bracket polynomial is even, and the span of the Jones polynomial belongs to~$\frac12\mathbb Z$.
\end{proof}

\begin{figure}[ht]
	\begin{center}
		\begin{tikzpicture}[scale=0.9]
			\pic[
			rotate=90,
			braid/.cd,
			every strand/.style={ultra thick},
			strand 1/.style={black},
			strand 2/.style={black}]
			at (0,1) {braid={s_1^{1} s_1^{1}}};
			\node[black] at (3.5, 2) {$\cdots$};
			\node[black] at (3.5, 1) {$\cdots$};
			\pic[
			rotate=90,
			braid/.cd,
			every strand/.style={ultra thick},
			strand 1/.style={black},
			strand 2/.style={black}]
			at (4.5,1) {braid={s_1^{1} }};
			\draw  [ultra thick, black]  (0,3) -- (7,3);
			\draw  [ultra thick, black]  (6,1) -- (6.5,1);
			\draw  [ultra thick, black]  (6.6,2.1) -- (7,2.1);
			\draw  [ultra thick, black]  (6.49,1) -- (6.49,2.5);
			\draw[ultra thick, black] (0,3) arc (90:270: 0.25 and 0.435);
			\draw[ultra thick, black] (7,3) arc (90:-90:0.25 and 0.435);
			\draw [black, very thick,->](0.1,1) -- (0.3,1);
			\fill[black] (0,1) circle (3pt);
			\fill[black] (6.49,2.5) circle (3pt);
			\node[black] at (0.5, 1.5) {$c_1$};
			\draw[black, ultra thick, decorate, decoration={brace, amplitude=6pt, mirror}]
			(0.3,0.8) -- (5.8,0.8) node [midway, yshift=-0.4cm] {$m-1$ crossings};
			\node[black, above] at (3.1,-1.0) {(a)};
		\end{tikzpicture}
		\quad
		\begin{tikzpicture}[scale=0.9]
			\pic[
			rotate=90,
			braid/.cd,
			every strand/.style={ultra thick},
			strand 1/.style={black},
			strand 2/.style={black}]
			at (0,1) {braid={s_1^{1} s_1^{1}}};
			\node[black] at (3.5, 2) {$\cdots$};
			\node[black] at (3.5, 1) {$\cdots$};
			\pic[
			rotate=90,
			braid/.cd,
			every strand/.style={ultra thick},
			strand 1/.style={black},
			strand 2/.style={black}]
			at (4.5,1) {braid={s_1^{1} }};
			\draw  [ultra thick, black]  (6,2.1) -- (6.5,2.1);
			\draw  [ultra thick, black]  (6,1) -- (6.7,1);
			\draw  [ultra thick, black]  (0,3) -- (6.7,3);
			\draw  [ultra thick, black]  (-0.4,1) -- (0,1);
			\draw[ultra thick, black] (6.7,3) arc (90:-90: 0.25 and 0.98);
			\draw[ultra thick, black] (0,3) arc (90:270: 0.25 and 0.7);
			\draw[ultra thick, black] (-0.4,2.1) arc (90:270:0.25 and 0.55);
			\draw[black, very thick,->] (6,2.1)-- (5.8,2.1);
			\fill[black] (0,1.57) circle (3pt);
			\fill[black] (6.5,2.1) circle (3pt);
			\draw[black, ultra thick, decorate, decoration={brace, amplitude=6pt, mirror}]
			(0.3,0.8) -- (6,0.8) node [midway, yshift=-0.4cm] {$m-1$ crossings};
			\node[black, above] at (3.1,-1.0) {(b)};
			\end{tikzpicture}
			\end{center}
			\caption{Multi-knotoids $\mathcal N(m)$ and $\mathcal L(m)$.} \label{fig11}
			\end{figure}
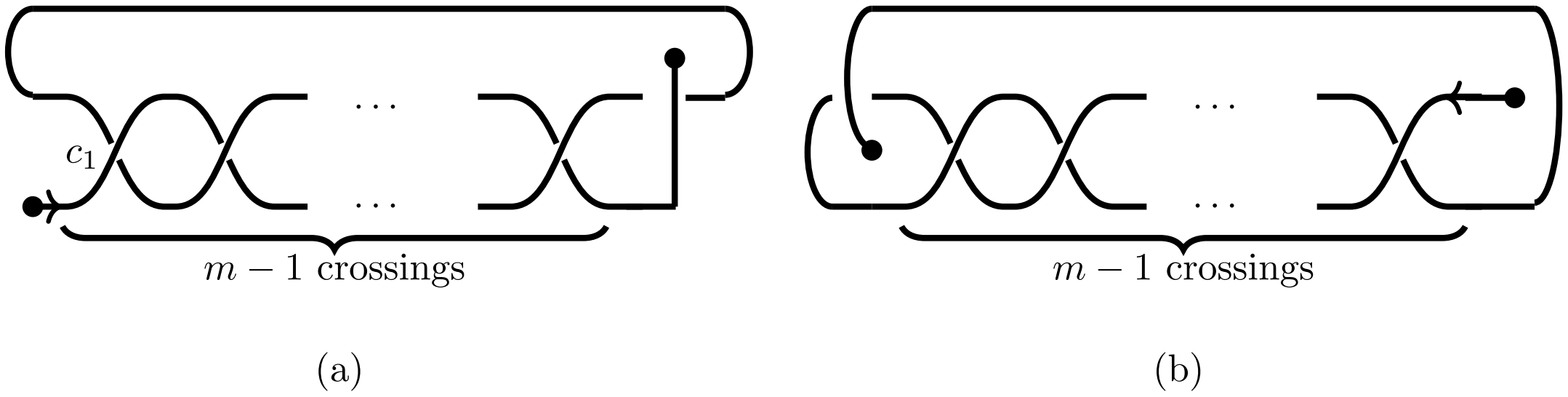

Consider two infinite families of spherical multi-knotoid diagrams: $\mathcal N(m)$, $m \geq 1$, see Fig.~\ref{fig11}(a), and $\mathcal L(m)$, $m \geq 1$, see Fig.~\ref{fig11}(b). 
			
			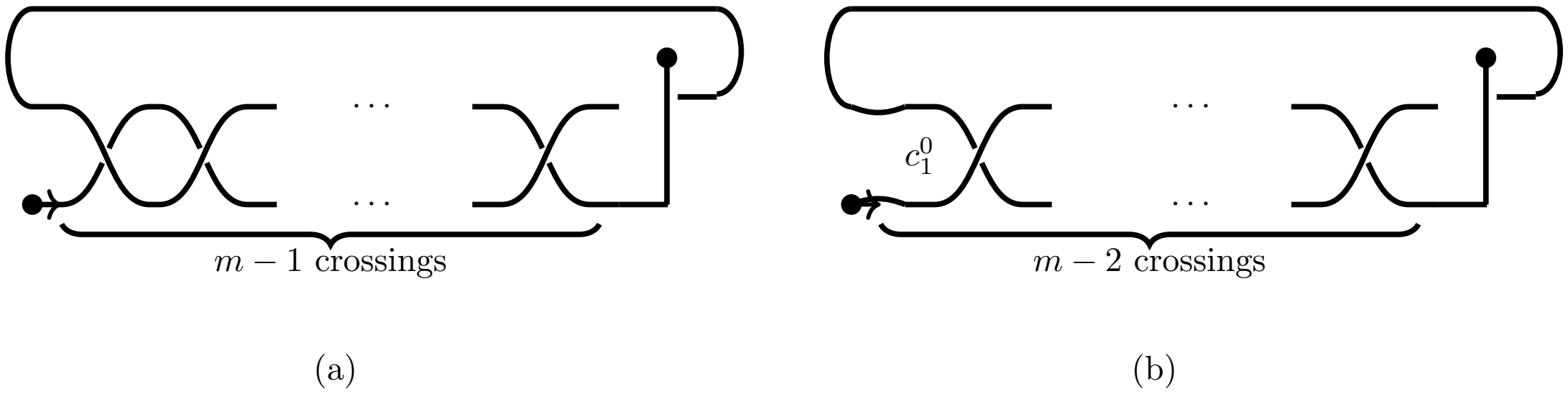
\begin{figure}[ht]
				\begin{center}
					\scalebox{0.9}	{
						\begin{tikzpicture}
							\pic[
							rotate=90,
							braid/.cd,
							every strand/.style={ultra thick},
							strand 1/.style={black},
							strand 2/.style={black}]
							at (0,1) {braid={s_1^{-1} s_1^{1}}};
							\node[black] at (3.5, 2) {$\cdots$};
							\node[black] at (3.5, 1) {$\cdots$};
							\pic[
							rotate=90,
							braid/.cd,
							every strand/.style={ultra thick},
							strand 1/.style={black},
							strand 2/.style={black}]
							at (4.5,1) {braid={s_1^{1} }};
							\draw  [ultra thick, black]  (0,3) -- (7,3);
							\draw  [ultra thick, black]  (6,1) -- (6.5,1);
							\draw  [ultra thick, black]  (6.6,2.1) -- (7,2.1);
							\draw  [ultra thick, black]  (6.49,1) -- (6.49,2.5);
							\draw[ultra thick, black] (0,3) arc (90:270: 0.25 and 0.5);
							\draw[ultra thick, black] (7,3) arc (90:-90:0.25 and 0.435);
							\draw [black, very thick,->](0.1,1) -- (0.3,1);
							\fill[black] (0,1) circle (3pt);
							\fill[black] (6.49,2.5) circle (3pt);
							\draw[black, ultra thick, decorate, decoration={brace, amplitude=6pt, mirror}]
							(0.3,0.8) -- (5.8,0.8) node [midway, yshift=-0.4cm] {$m-1$ crossings};
							\node[black, above] at (3.1,-1.0) {(a)};
						\end{tikzpicture}
						\qquad
						\begin{tikzpicture}
							\draw[ultra thick, black] (0,1) .. controls (0.2,1.08) and (0.35,1.08) .. (0.55,1);
							\draw[ultra thick, black] (0,2) .. controls (0.2,1.92) and (0.35,1.92) .. (0.55,2);
							\pic[
							rotate=90,
							braid/.cd,
							every strand/.style={ultra thick},
							strand 1/.style={black},
							strand 2/.style={black}
							] at (0.55,1) {braid={s_1^{1}}};
							\node[black] at (3.5, 2) {$\cdots$};
							\node[black] at (3.5, 1) {$\cdots$};
							\pic[
							rotate=90,
							braid/.cd,
							every strand/.style={ultra thick},
							strand 1/.style={black},
							strand 2/.style={black}
							] at (4.5,1) {braid={s_1^{1}}};		
							\draw [ultra thick, black] (0,3) -- (7,3);
							\draw [ultra thick, black] (6,1) -- (6.5,1);
							\draw [ultra thick, black] (6.6,2.1) -- (7,2.1);
							\draw [ultra thick, black] (6.49,1) -- (6.49,2.5);
							\draw[ultra thick, black] (0,3) arc (90:270: 0.25 and 0.5);
							\draw[ultra thick, black] (7,3) arc (90:-90:0.25 and 0.435);
							\draw [black, very thick,->](0.1,1) -- (0.3,1);
							\fill[black] (0,1) circle (3pt);
							\fill[black] (6.49,2.5) circle (3pt);
							\node[black] at (0.7, 1.5) {$c_1^0$};
							\draw[black, ultra thick, decorate, decoration={brace, amplitude=6pt, mirror}]
							(0.3,0.8) -- (5.8,0.8) node [midway, yshift=-0.4cm] {$m-2$ crossings};
							\node[black, above] at (3.1,-1.0) {(b)};
						\end{tikzpicture}
					}
				\end{center}		
				\caption{$K_+ = \mathcal N(m-2)$ and $K_0 = \mathcal N (m-1)$ from the Conway triple for $K_- = \mathcal N(m)$.}\label{fig16}
			\end{figure}
			
			Note that the multi-knotoid $\mathcal N(2p)$, $p > 1$, is a knotoid. 	
			\begin{lemma} \label{lemma3.8}
				For the spherical knotoid $\mathcal N(2p)$, $p > 1$, we have $\mathbf{span}(V_{\mathcal N(2p)}(t)) = p$.
			\end{lemma}
			
			\begin{proof}
				As seen from Fig.~\ref{fig11}(a), the signs of all crossings are $-1$. We use the skein relation~(\ref{eqn6}). Consider $c_1$, the first crossing encountered when moving from the head of the knotoid. Setting $K_- = \mathcal N(m)$ in the notation of the Conway triple, we obtain $K_+ = \mathcal N(m-2)$ (after applying a type II Reidemeister move) and $K_0 = \mathcal N(m-1)$, shown in Fig.~\ref{fig16}. For the crossing $c_1^0$ of the diagram $K_0$, we write relation~(\ref{eqn8}). Then
				\begin{equation}\label{eqn9}
					\langle K_0 \rangle_\circ = \langle \mathcal N(m-1) \rangle_\circ = -A^{4}\langle U \rangle_\circ + (-A^{2})\langle \mathcal N(m-2) \rangle_\circ, 
				\end{equation}
				where $U$ denotes the trivial knotoid. Denote $x_{m} = \langle \mathcal N(m) \rangle_\circ$, where $x_0 = \langle U  \rangle_\circ = 1$. Then, applying relation~(\ref{eqn6}) to the crossing $c_1$ of the knotoid $\mathcal N(m)$, we obtain 
				\begin{equation}\label{eqn13}
					-A^4\cdot x_{m-2} + A^{-4}\cdot x_{m} =  (A^2-A^{-2})\cdot x_{m-1}.
				\end{equation}
				Substitute $x_{m-1} = - A^4 - A^2 x_{m-2}$, obtained from~(\ref{eqn9}), into~(\ref{eqn13}) and simplify. Then
				\begin{equation}\label{eqn10}
					x_{m} = A^4\cdot x_{m-2} - A^{10} + A^6.
				\end{equation}
				
				Note that $x_2 = A^4 - A^{10} + A^6$. Further, for all $x_m$ with $m > 2$, the minimum degree of $A$ is $6$, while the maximum degree of $A$ increases by $4$ when passing from $x_m$ to $x_{m+2}$ for $m > 1$. Since the number of steps from $x_2$ to $x_m$ is ${{m-2}\over{2}}$, the maximum degree of $A$ for $x_m$ is $2(m-2)+10 = 2m +6$. Since $m=2p$, we obtain $\mathbf{span}(V_{\mathcal N(2p)}(t)) = {{4p+6-6}\over{4}} = p$. 
			\end{proof}
			
			Note that the multi-knotoid $\mathcal L(m)$, $m \geq 2$, is a knotoid. 	
			\begin{lemma} \label{lemma3.9}
				For the spherical knotoid $\mathcal L(m)$, $m \geq 2$, we have $\mathbf{span}(V_{\mathcal L(m)} (t)) = m - \frac{1}{2}$. 
			\end{lemma}
			
			\begin{proof} 
				Consider the odd case $m=2p+1$, $p \geq 1$. We use the skein relations. Choose the crossing enclosed by the dashed circle in Fig.~\ref{fig15}(a) and consider the Conway triple for $K_+ = \mathcal L(2p+1)$. Then $K_- = \mathcal L(2p-1)$, and $K_0$ is shown in Fig.~\ref{fig15}(b). It is easy to see that $\langle K_0 \rangle_\circ = -A^4-A^2$.
				\begin{figure}[ht]
					\begin{center}
						\begin{tikzpicture}[scale=0.9, baseline=(current bounding box.center)]
							\pic[
							rotate=90,
							braid/.cd,
							every strand/.style={ultra thick},
							strand 1/.style={black},
							strand 2/.style={black}]
							at (0,1) {braid={s_1^{1} s_1^{1}}};
							\node[black] at (3.5, 2) {$\cdots$};
							\node[black] at (3.5, 1) {$\cdots$};
							\pic[
							rotate=90,
							braid/.cd,
							every strand/.style={ultra thick},
							strand 1/.style={black},
							strand 2/.style={black}]
							at (4.5,1) {braid={s_1^{1} }};
							\draw [ultra thick, black] (6,2.1) -- (6.5,2.1);
							\draw [ultra thick, black] (6,1)   -- (6.7,1);
							\draw [ultra thick, black] (0,3)   -- (6.7,3);
							\draw [ultra thick, black] (-0.4,1) -- (0,1);
							\draw[ultra thick, black] (6.7,3) arc (90:-90: 0.25 and 0.98);
							\draw[ultra thick, black] (0,3)   arc (90:270: 0.25 and 0.7);
							\draw[ultra thick, black] (-0.4,2.1) arc (90:270: 0.25 and 0.55);
							\draw [black, very thick,->] (6,2.1) -- (5.8,2.1);
							\fill[black] (0,1.57)  circle (3pt);
							\fill[black] (6.5,2.1) circle (3pt);
							\draw[black, ultra thick, decorate, decoration={brace, amplitude=6pt, mirror}]
							(0.3,0.8) -- (6,0.8) node [midway, yshift=-0.4cm] {$m-1$ crossings};
							\draw[ultra thick, red, dashed] (0.8,1.5) circle[x radius=0.5, y radius=0.5];
							\node[black, above] at (3.1,-1.0) {(a)};
						\end{tikzpicture}
						\hspace{0.5cm}
						\begin{tikzpicture}[scale=0.8, baseline=(current bounding box.center)]
							\draw [ultra thick, black] (5,3)   -- (7,3);
							\draw [ultra thick, black] (5,2)   -- (6.4,2);
							\draw [ultra thick, black] (6.6,2) -- (7,2);
							\draw [ultra thick, black] (6.49,1) -- (6.49,2.5);
							\draw[ultra thick, black] (5,3) arc (90:270: 0.25 and 0.5);
							\draw[ultra thick, black] (7,3) arc (90:-90: 0.25 and 0.5);
							\draw [black, very thick,->] (6.49,1) -- (6.49,1.5);
							\draw [black, very thick,->] (7,3) -- (6,3);
							\fill[black] (6.49,1)   circle (3pt);
							\fill[black] (6.49,2.5) circle (3pt);
							\node[black, above] at (6.1,-1.0) {(b)};
						\end{tikzpicture}
						\hspace{0.5cm}
						\begin{tikzpicture}[scale=0.8, baseline=(current bounding box.center)]
							\draw [ultra thick, black] (5,3)   -- (7,3);
							\draw [ultra thick, black] (5,2)   -- (6.4,2);
							\draw [ultra thick, black] (6.6,2) -- (7,2);
							\draw [ultra thick, black] (6.49,1) -- (6.49,2.5);
							\draw[ultra thick, black] (5,3) arc (90:270: 0.25 and 0.5);
							\draw[ultra thick, black] (7,3) arc (90:-90: 0.25 and 0.5);
							\draw [black, very thick,->] (6.49,1) -- (6.49,1.5);
							\draw [black, very thick,->] (5,3) -- (6,3);
							\fill[black] (6.49,1)   circle (3pt);
							\fill[black] (6.49,2.5) circle (3pt);
							\node[black, above] at (6.1,-1.0) {(c)};
						\end{tikzpicture}
					\end{center}
					\caption{The knotoid $\mathcal L(m)$ and the multi-knotoid $K_0$ for the even and odd cases.} \label{fig15}
				\end{figure}
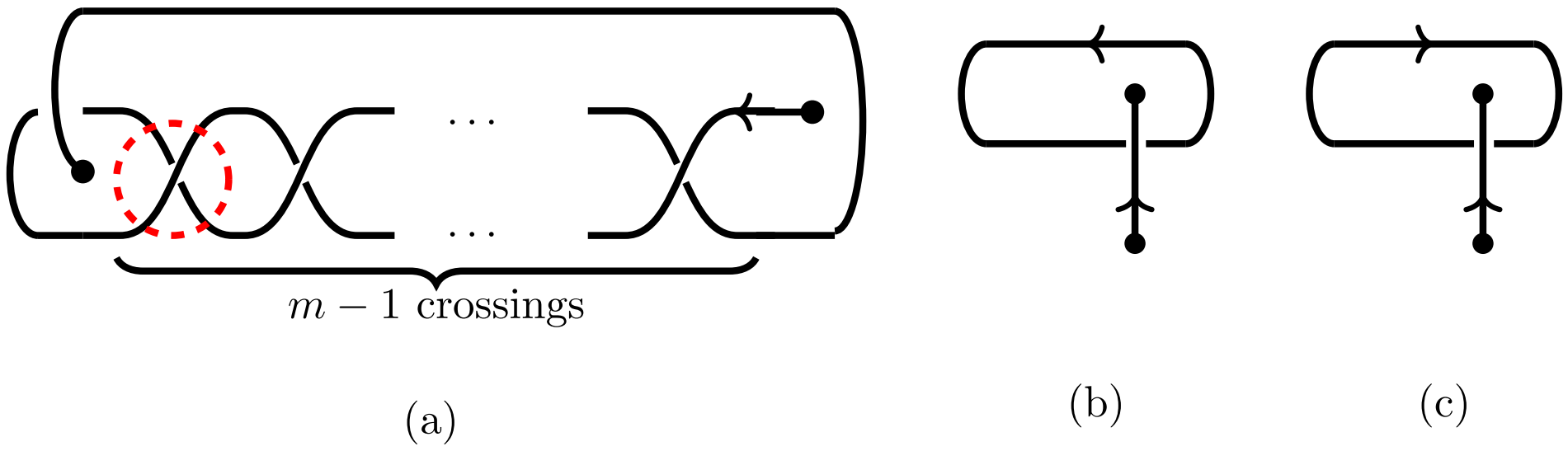
				
				Denote $y_{2p+1} = \langle \mathcal L(2p+1) \rangle_\circ$, where $y_1 = \langle U \rangle_\circ = 1$. Then, using the skein relation~(\ref{eqn6}) and the value of $\langle K_0 \rangle_\circ$, we obtain the recurrence relation:  
				\begin{equation}\label{eqn11}
					y_{2p+1} = A^{-8}y_{2p-1}+A^2+1-A^{-2}-A^{-4}. 
				\end{equation}
				Since the degrees of $A$ appearing in $y_{2p-1}$ decrease by $8$ when passing to $y_{2p+1}$ in relation~(\ref{eqn11}), the maximum degree remains equal to $2$. It is also easy to see that the minimum degree of $A$ is $-8p$. Hence,
				$$
				\mathbf{span}(V_{ \mathcal L(2p+1)}(t)) = {{2+8p}\over{4}} = 2p+\frac{1}{2} = m-1+\frac{1}{2} = m-\frac{1}{2}.
				$$ 
				
				Consider the even case $m=2p+2$, $p \geq 0$. Again, choose the crossing enclosed by the dashed circle in Fig.~\ref{fig15}(a) and consider the Conway triple for $K_+ = \mathcal L(2p+2)$. Then $K_- = \mathcal L(2p)$, and $K_0$ is shown in Fig.~\ref{fig15}(c). It is easy to see that $\langle K_0 \rangle_\circ = -A^{-4}-A^{-2}$.
				
				Denote $y_{2p+2} = \langle \mathcal L(2p+2) \rangle_\circ$, where $y_0 = \langle U \rangle_\circ = 1$. Then, using the skein relation~(\ref{eqn6}) and the value of $\langle K_0 \rangle_\circ$, we obtain the recurrence relation:  
				\begin{equation}\label{eqn12}
					y_{2p+2} = A^{-8}y_{2p}+A^{-4}+A^{-6}-A^{-8}-A^{-10}. 
				\end{equation}
				Note that $y_2 = -A^{-10}+A^{-6}+A^{-4}$. Further, for all $y_{2p}$ with $p>0$, the maximum degree of $A$ is $-4$, while the minimum degree of $A$ decreases by $8$ when passing from $y_{2p}$ to $y_{2p+2}$ for $p>0$. Since the number of steps from $y_2$ to $y_{2p+2}$ is $p$, the minimum degree of $A$ for $y_{2p+2}$ is $-8p-10$. Hence,
				$$
				\mathbf{span}(V_{ \mathcal L(2p+2)}(t)) = {{6+8p}\over{4}} = 2p+1+\frac{1}{2}.
				$$  
				
				Recall that in this case $m=2p+2$. Thus, for any integer $m\geq2$, we have $\mathbf{span}(V_{\mathcal L(m)} (t)) = m-\frac{1}{2}$.
				\end{proof}
				
				\begin{lemma} \label{lemma3.10}
					Let $\mathcal L(m)$, $m \geq 2$, be the spherical knotoid shown in Fig.~\ref{fig11}(b). If $m$ is odd, then $\mathcal L(m)$ can be transformed into the knotoid $\mathcal N(m-1)$ by a single $H(2)$-move, whereas if $m$ is even, it can be transformed into the knotoid $\phi(T(m-1,2))$. Here $T(m-1,2)$ is the torus knot in Fig.~\ref{fig2}, and $\phi$ is Turaev's map from the set of knots to the set of spherical knotoids. 
				\end{lemma}
				
				\begin{proof} 
					In the odd case $m = 2p+1$, we apply the $H(2)$-move shown in Fig.~\ref{fig12}.
					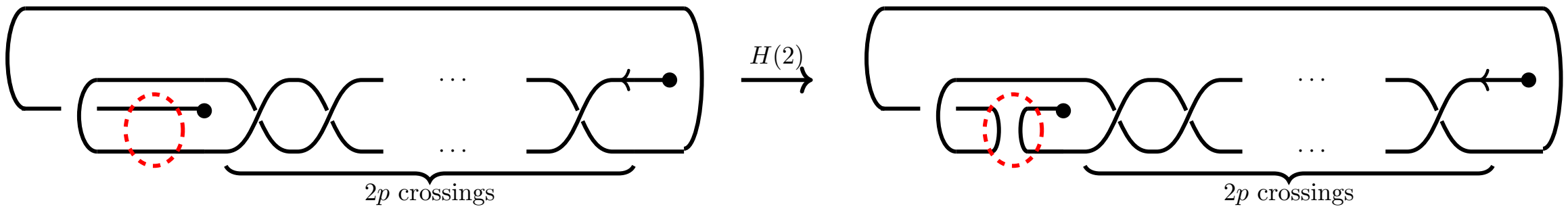
\begin{figure}[ht]
						\begin{center} 
							\scalebox{0.7}	{
								\begin{tikzpicture} [ baseline=(current bounding box.center)]
									\begin{scope}
										\pic[
										rotate=90,
										braid/.cd,
										every strand/.style={ultra thick},
										strand 1/.style={black},
										strand 2/.style={black}]
										at (0,1) {braid={s_1^{1} s_1^{1}}};
										\node[black] at (3.5, 2) {$\cdots$};
										\node[black] at (3.5, 1) {$\cdots$};
										\pic[
										rotate=90,
										braid/.cd,
										every strand/.style={ultra thick},
										strand 1/.style={black},
										strand 2/.style={black}]
										at (4.5,1) {braid={s_1^{1} }};
										\draw [ultra thick, black] (6,2) -- (6.5,2);
										\draw [ultra thick, black] (-1.5,2)   -- (0,2);
										\draw [ultra thick, black] (-1.5,1)   -- (0,1);
										\draw [ultra thick, black] (-1.5,1.6)   -- (0,1.6);
										\draw [ultra thick, black] (-2.55,1.6)   -- (-2,1.6);
										\draw [ultra thick, black] (6,1)   -- (6.7,1);
										\draw [ultra thick, black] (-2.5,3)   -- (6.7,3);
										\draw [ultra thick, black] (-0.4,1) -- (0,1);
										\draw[ultra thick, black] (6.7,3) arc (90:-90: 0.25 and 0.98);
										\draw[ultra thick, black] (-2.5,3)   arc (90:270: 0.25 and 0.7);
										\draw[ultra thick, black] (-1.5,2) arc (90:270: 0.25 and 0.5);
										\draw [black, very thick,->] (6,2) -- (5.8,2);
										\fill[black] (0,1.57)  circle (3pt);
										\fill[black] (6.5,2) circle (3pt);
										\draw[black, ultra thick, decorate, decoration={brace, amplitude=6pt, mirror}]
										(0.3,0.8) -- (6,0.8) node [midway, yshift=-0.4cm] {$2p$ crossings};
										\draw[ultra thick, red, dashed] (-0.7,1.3) circle[x radius=0.4, y radius=0.5];
									\end{scope}
									\draw[->, ultra thick] (7.5, 2) -- (8.5, 2) node[midway, above] {$H(2)$};
									\begin{scope}[xshift=12cm]
										\pic[
										rotate=90,
										braid/.cd,
										every strand/.style={ultra thick},
										strand 1/.style={black},
										strand 2/.style={black}]
										at (0,1) {braid={s_1^{1} s_1^{1}}};
										\node[black] at (3.5, 2) {$\cdots$};
										\node[black] at (3.5, 1) {$\cdots$};
										\pic[
										rotate=90,
										braid/.cd,
										every strand/.style={ultra thick},
										strand 1/.style={black},
										strand 2/.style={black}]
										at (4.5,1) {braid={s_1^{1} }};
										\draw [ultra thick, black] (6,2) -- (6.5,2);
										\draw [ultra thick, black] (-1.5,2)   -- (0,2);
										\draw [ultra thick, black] (-1.5,1)   -- (-1,1);
										\draw [ultra thick, black] (-0.5,1)   -- (0,1);
										\draw [ultra thick, black] (-1.5,1.6)   -- (-1,1.6);
										\draw [ultra thick, black] (-0.5,1.6)   -- (0,1.6);
										\draw [ultra thick, black] (-2.55,1.6)   -- (-2,1.6);
										\draw [ultra thick, black] (6,1)   -- (6.7,1);
										\draw [ultra thick, black] (-2.5,3)   -- (6.7,3);
										\draw [ultra thick, black] (-0.4,1) -- (0,1);
										\draw[ultra thick, black] (6.7,3) arc (90:-90: 0.25 and 0.98);
										\draw[ultra thick, black] (-2.5,3)   arc (90:270: 0.25 and 0.7);
										\draw[ultra thick, black] (-1.5,2) arc (90:270: 0.25 and 0.5);
										\draw[ultra thick, black] (-1,1.6) arc (90:-90: 0.1 and 0.3);
										\draw[ultra thick, black] (-0.5,1.6) arc (90:270: 0.1 and 0.3);
										\draw [black, very thick,->] (6,2) -- (5.8,2);
										\fill[black] (0,1.57)  circle (3pt);
										\fill[black] (6.5,2) circle (3pt);
										\draw[black, ultra thick, decorate, decoration={brace, amplitude=6pt, mirror}]
										(0.3,0.8) -- (6,0.8) node [midway, yshift=-0.4cm] {$2p$ crossings};
										\draw[ultra thick, red, dashed] (-0.7,1.3) circle[x radius=0.4, y radius=0.5];
									\end{scope}
								\end{tikzpicture}
							}
						\end{center}
						\caption{An $H(2)$-move transforming $\mathcal L(2p+1)$ into $\mathcal N(2p)$.} \label{fig12}
					\end{figure}
					
					\noindent 		
					In the even case $m = 2p+2$, we apply the $H(2)$-move shown in Fig.~\ref{fig13}.
					\begin{figure}[ht]
						\begin{center} 
							\scalebox{0.7}	{
								\begin{tikzpicture} [ baseline=(current bounding box.center)]
									\begin{scope}
										\pic[
										rotate=90,
										braid/.cd,
										every strand/.style={ultra thick},
										strand 1/.style={black},
										strand 2/.style={black}]
										at (0,1) {braid={s_1^{1} s_1^{1}}};
										\node[black] at (3.5, 2) {$\cdots$};
										\node[black] at (3.5, 1) {$\cdots$};
										\pic[
										rotate=90,
										braid/.cd,
										every strand/.style={ultra thick},
										strand 1/.style={black},
										strand 2/.style={black}]
										at (4.5,1) {braid={s_1^{1} }};
										\draw [ultra thick, black] (6,2) -- (6.5,2);
										\draw [ultra thick, black] (-1.5,2)   -- (0,2);
										\draw [ultra thick, black] (-1.5,1)   -- (0,1);
										\draw [ultra thick, black] (-1.5,1.6)   -- (0,1.6);
										\draw [ultra thick, black] (-2.55,1.6)   -- (-2,1.6);
										\draw [ultra thick, black] (6,1)   -- (6.7,1);
										\draw [ultra thick, black] (-2.5,3)   -- (6.7,3);
										\draw [ultra thick, black] (-0.4,1) -- (0,1);
										\draw[ultra thick, black] (6.7,3) arc (90:-90: 0.25 and 0.98);
										\draw[ultra thick, black] (-2.5,3)   arc (90:270: 0.25 and 0.7);
										\draw[ultra thick, black] (-1.5,2) arc (90:270: 0.25 and 0.5);
										\draw [black, very thick,->] (6,2) -- (5.8,2);
										\fill[black] (0,1.57)  circle (3pt);
										\fill[black] (6.5,2) circle (3pt);
										\draw[black, ultra thick, decorate, decoration={brace, amplitude=6pt, mirror}]
										(0.3,0.8) -- (6,0.8) node [midway, yshift=-0.4cm] {$2p+1$ crossings};
										\draw[ultra thick, red, dashed] (-0.7,1.7) circle[x radius=0.4, y radius=0.5];
									\end{scope}
									\draw[->, ultra thick] (7.5, 2) -- (8.5, 2) node[midway, above] {$H(2)$}; 
									\begin{scope}[xshift=12cm]
										\pic[
										rotate=90,
										braid/.cd,
										every strand/.style={ultra thick},
										strand 1/.style={black},
										strand 2/.style={black}]
										at (0,1) {braid={s_1^{1} s_1^{1}}};
										\node[black] at (3.5, 2) {$\cdots$};
										\node[black] at (3.5, 1) {$\cdots$};
										\pic[
										rotate=90,
										braid/.cd,
										every strand/.style={ultra thick},
										strand 1/.style={black},
										strand 2/.style={black}]
										at (4.5,1) {braid={s_1^{1} }};
										\draw [ultra thick, black] (6,2) -- (6.5,2);
										\draw [ultra thick, black] (-1.5,2)   -- (-1,2);
										\draw [ultra thick, black] (-0.5,2)   -- (0,2);
										\draw [ultra thick, black] (-1.5,1)   -- (0,1);
										\draw [ultra thick, black] (-1.5,1.6)   -- (-1,1.6);
										\draw [ultra thick, black] (-0.5,1.6)   -- (0,1.6);
										\draw [ultra thick, black] (-2.55,1.6)   -- (-2,1.6);
										\draw [ultra thick, black] (6,1)   -- (6.7,1);
										\draw [ultra thick, black] (-2.5,3)   -- (6.7,3);
										\draw [ultra thick, black] (-0.4,1) -- (0,1);
										\draw[ultra thick, black] (6.7,3) arc (90:-90: 0.25 and 0.98);
										\draw[ultra thick, black] (-2.5,3)   arc (90:270: 0.25 and 0.7);
										\draw[ultra thick, black] (-1.5,2) arc (90:270: 0.25 and 0.5);
										\draw[ultra thick, black] (-1,2) arc (90:-90: 0.2 and 0.2);
										\draw[ultra thick, black] (-0.5,2) arc (90:270: 0.2 and 0.2);
										\draw [black, very thick,->] (6,2) -- (5.8,2);
										\fill[black] (0,1.57)  circle (3pt);
										\fill[black] (6.5,2) circle (3pt);
										\draw[black, ultra thick, decorate, decoration={brace, amplitude=6pt, mirror}]
										(0.3,0.8) -- (6,0.8) node [midway, yshift=-0.4cm] {$2p+1$ crossings};
										\draw[ultra thick, red, dashed] (-0.7,1.8) circle[x radius=0.4, y radius=0.5];
									\end{scope}
								\end{tikzpicture}
							}
						\end{center}
						\caption{An $H(2)$-move transforming $\mathcal L(2p+2)$ into $\phi(T(2p+1,2))$.} \label{fig13}
					\end{figure}
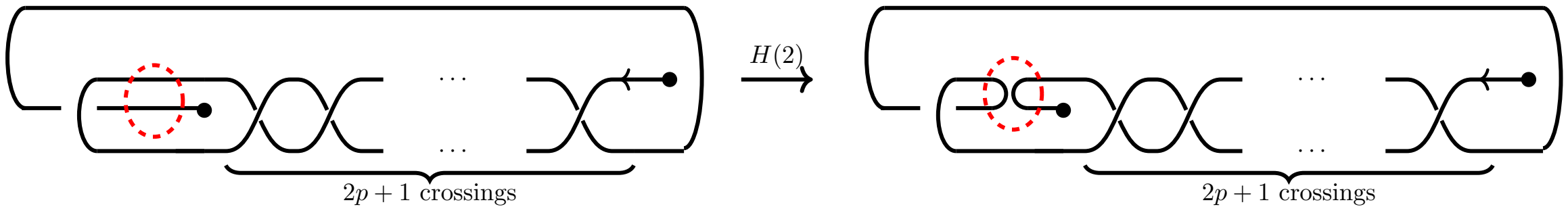
					\noindent 
					The assertion of the lemma follows from these moves.
				\end{proof}
				
				In~\cite{FMS}, the following lemma was proved.
				\begin{lemma}~\cite[Cor.~3.4]{FMS} \label{lemma3.14}
					Let $T(2p+1,2)$ be the torus knot shown in Fig.~\ref{fig2}. Then $\mathbf{span}(V_{T(2p+1,2)} (t)) = 2p+1$. 
				\end{lemma}
				
				\begin{figure}[ht]
					\begin{center} 
							\begin{tikzpicture} 
								\pic[
								rotate=90,
								braid/.cd,
								every strand/.style={ultra thick},
								strand 1/.style={black},
								strand 2/.style={black}] 
								at (0,1) {braid={s_1^{1} s_1^{1} }}; 
								\node[black] at (3.5, 2) {$\cdots$};
								\node[black] at (3.5, 1) {$\cdots$};
								\pic[
								rotate=90,
								braid/.cd,
								every strand/.style={ultra thick},
								strand 1/.style={black},
								strand 2/.style={black}] 
								at (4.5,1) {braid={s_1^{1} s_1^{1} }}; 
								\draw  [ultra thick, black]  (0,0) -- (7,0);
								\draw  [ultra thick, black]  (0,3) -- (7,3);
								\draw[ultra thick, black] (0,1) arc (90:270: 0.25 and 0.5);
								\draw[ultra thick, black] (0,3) arc (90:270: 0.25 and 0.5);
								\draw[ultra thick, black] (7,1) arc (90:-90:0.25 and 0.5);
								\draw[ultra thick, black] (7,3) arc (90:-90:0.25 and 0.5);
								\draw [ultra thick, decorate, decoration={brace, amplitude=6pt, mirror}] 
								(0.5,0.8) -- (6.5,0.8) node [midway, yshift=-0.5cm] {$2p+1$ crossings};
							\end{tikzpicture}
					\end{center}
					\caption{The torus knot $T(2p+1,2)$.} \label{fig2}
				\end{figure}
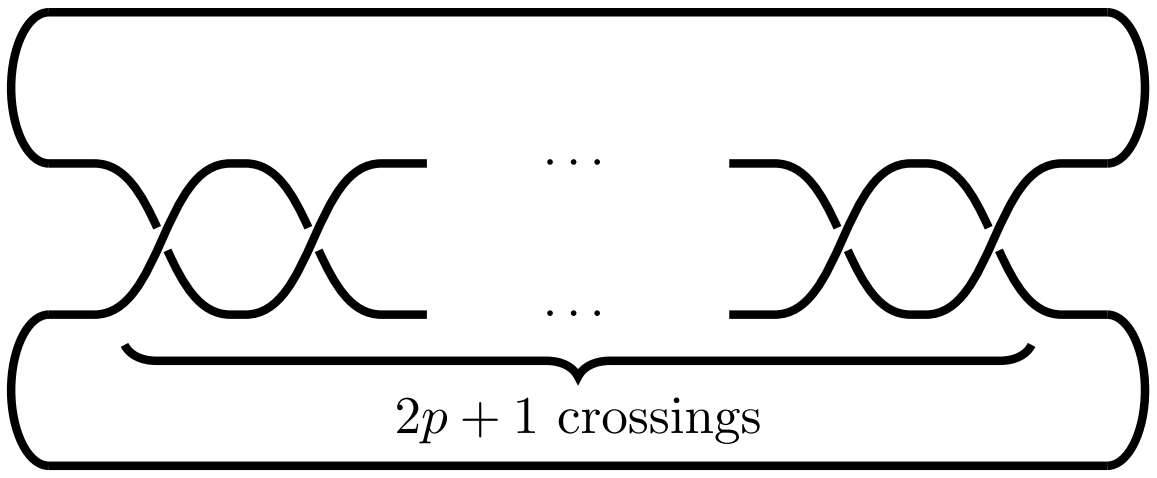
				
				\begin{theorem}\label{theorem1.6}
					Denote by $G(H(n), \mathbf{span}(V))$ the quotient graph of the $H(n)$-Gordian graph of knotoids by the span $\mathbf{span}(V(t))$ of the Jones polynomial $V(t)$. Then the graph $G(H(n), \mathbf{span}(V))$ is Gromov hyperbolic.
				\end{theorem}
				
				\begin{proof} 
					It is easy to see that the knotoids $\phi(T(2p+1,2))$ and $\mathcal N(2p)$ can be transformed into the trivial knotoid by a single $H(2)$-move (see Fig.~\ref{fig17} and Fig.~\ref{fig18}, respectively). For $p \in \frac{1}{2} \mathbb Z$, we introduce the class of knotoids whose Jones polynomial span is equal to $p$:
					$$
					[K]_p = \{ K : \mathbf{span}(V_K(t)) = p\}. 
					$$
					
					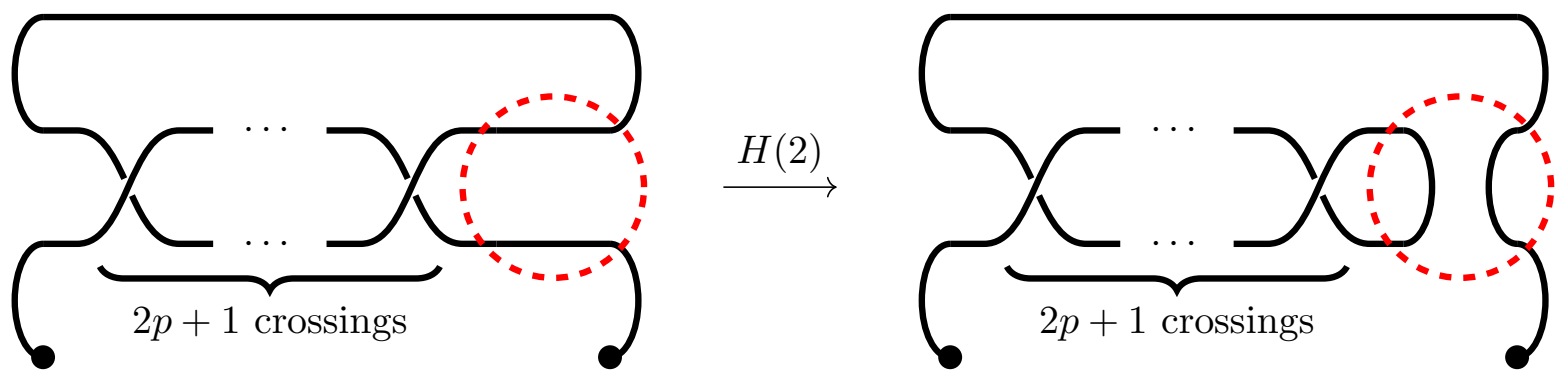
\begin{figure}[ht]
						\begin{center} 
							\scalebox{0.75}	{
								\begin{tikzpicture} [ baseline=(current bounding box.center)]
									\pic[
									rotate=90,
									braid/.cd,
									every strand/.style={ultra thick},
									strand 1/.style={black},
									strand 2/.style={black}] 
									at (0,1) {braid={s_1  }}; 
									\node[black] at (2, 2) {$\cdots$};
									\node[black] at (2, 1) {$\cdots$};
									\pic[
									rotate=90,
									braid/.cd,
									every strand/.style={ultra thick},
									strand 1/.style={black},
									strand 2/.style={black}] 
									at (2.5,1) {braid={s_1  }}; 
									\fill[black] (0,0)  circle (3pt);
									\fill[black] (5,0)  circle (3pt);
									\draw  [ultra thick, black]  (0,3) -- (5,3);
									\draw  [ultra thick, black]  (4,1) -- (5,1);
									\draw  [ultra thick, black]  (4,2) -- (5,2);
									\draw[ultra thick, black] (0,1) arc (90:270: 0.25 and 0.5);
									\draw[ultra thick, black] (0,3) arc (90:270: 0.25 and 0.5);
									\draw[ultra thick, black] (5,1) arc (90:-90:0.25 and 0.5);
									\draw[ultra thick, black] (5,3) arc (90:-90:0.25 and 0.5);
									\draw[ultra thick, red, dashed]  (4.5,1.5) circle (0.8);
									\draw [ultra thick, decorate, decoration={brace, amplitude=6pt, mirror}] 
									(0.5,0.8) -- (3.5,0.8) node [midway, yshift=-0.5cm] {$2p+1$ crossings};
									\draw[->] (6., 1.5) -- (7., 1.5) node[midway, above] {$H(2)$};
									\pic[
									rotate=90,
									braid/.cd,
									every strand/.style={ultra thick},
									strand 1/.style={black},
									strand 2/.style={black}] 
									at (8,1) {braid={s_1  }}; 
									\node[black] at (10, 2) {$\cdots$};
									\node[black] at (10, 1) {$\cdots$};
									\pic[
									rotate=90,
									braid/.cd,
									every strand/.style={ultra thick},
									strand 1/.style={black},
									strand 2/.style={black}] 
									at (10.5,1) {braid={s_1  }}; 
									\fill[black] (8,0)  circle (3pt);
									\fill[black] (13,0)  circle (3pt);
									\draw  [ultra thick, black]  (8,3) -- (13,3);
									\draw[ultra thick, black] (12,2) arc (90:-90: 0.25 and 0.5);
									\draw[ultra thick, black] (13,2) arc (90:270: 0.25 and 0.5);
									\draw[ultra thick, black] (8,1) arc (90:270: 0.25 and 0.5);
									\draw[ultra thick, black] (8,3) arc (90:270: 0.25 and 0.5);
									\draw[ultra thick, black] (13,1) arc (90:-90:0.25 and 0.5);
									\draw[ultra thick, black] (13,3) arc (90:-90:0.25 and 0.5);
									\draw[ultra thick, red, dashed]  (12.5,1.5) circle (0.8);
									\draw [ultra thick, decorate, decoration={brace, amplitude=6pt, mirror}] 
									(8.5,0.8) -- (11.5,0.8) node [midway, yshift=-0.5cm] {$2p+1$ crossings};
								\end{tikzpicture}
							}
						\end{center}
						\caption{An $H(2)$-move transforming $\phi(T(2p+1,2))$ into $U$.} \label{fig17}
					\end{figure}

					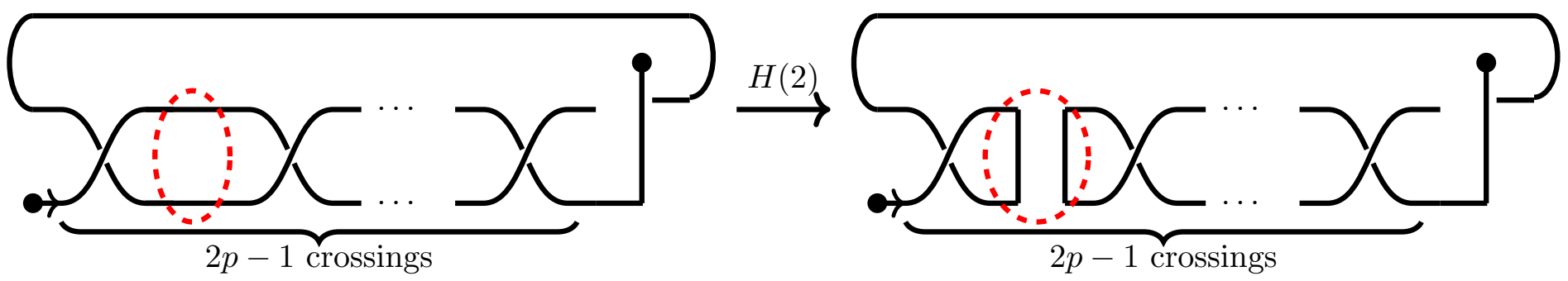
\begin{figure}[ht]
						\begin{center} 
							\scalebox{0.75}	{
								\begin{tikzpicture} [ baseline=(current bounding box.center)]
									\begin{scope}
										\pic[
										rotate=90,
										braid/.cd,
										every strand/.style={ultra thick},
										strand 1/.style={black},
										strand 2/.style={black}]
										at (0,1) {braid={s_1^{1}}};
										
										\draw  [ultra thick, black]  (1.2,1) -- (2,1);
										\draw  [ultra thick, black]  (1.2,2) -- (2,2);
										
										\pic[
										rotate=90,
										braid/.cd,
										every strand/.style={ultra thick},
										strand 1/.style={black},
										strand 2/.style={black}]
										at (2,1) {braid={s_1^{1}}};
										\node[black] at (3.9, 2) {$\cdots$};
										\node[black] at (3.9, 1) {$\cdots$};
										\pic[
										rotate=90,
										braid/.cd,
										every strand/.style={ultra thick},
										strand 1/.style={black},
										strand 2/.style={black}]
										at (4.5,1) {braid={s_1^{1} }};
										\draw  [ultra thick, black]  (0,3) -- (7,3);
										\draw  [ultra thick, black]  (6,1) -- (6.5,1);
										\draw  [ultra thick, black]  (6.6,2.1) -- (7,2.1);
										\draw  [ultra thick, black]  (6.49,1) -- (6.49,2.5);
										\draw[ultra thick, black] (0,3) arc (90:270: 0.25 and 0.5);
										\draw[ultra thick, black] (7,3) arc (90:-90:0.25 and 0.435);
										\draw [black, very thick,->](0.1,1) -- (0.3,1);
										\fill[black] (0,1) circle (3pt);
										\fill[black] (6.49,2.5) circle (3pt);
										\draw[black, ultra thick, decorate, decoration={brace, amplitude=6pt, mirror}]
										(0.3,0.8) -- (5.8,0.8) node [midway, yshift=-0.4cm] {$2p-1$ crossings};
										
										\draw[ultra thick, red, dashed]  (1.7,1.5) circle[x radius=0.4, y radius=0.7];
									\end{scope}
									\draw[->, ultra thick] (7.5, 2) -- (8.5, 2) node[midway, above] {$H(2)$};
									\begin{scope}[xshift=9cm]
										\pic[
										rotate=90,
										braid/.cd,
										every strand/.style={ultra thick},
										strand 1/.style={black},
										strand 2/.style={black}]
										at (0,1) {braid={s_1^{1}}};
										
										\draw  [ultra thick, black]  (1.5,1) -- (1.5,2);
										\draw  [ultra thick, black]  (2,1) -- (2,2);
										
										\pic[
										rotate=90,
										braid/.cd,
										every strand/.style={ultra thick},
										strand 1/.style={black},
										strand 2/.style={black}]
										at (2,1) {braid={s_1^{1}}};
										\node[black] at (3.9, 2) {$\cdots$};
										\node[black] at (3.9, 1) {$\cdots$};
										\pic[
										rotate=90,
										braid/.cd,
										every strand/.style={ultra thick},
										strand 1/.style={black},
										strand 2/.style={black}]
										at (4.5,1) {braid={s_1^{1} }};
										\draw  [ultra thick, black]  (0,3) -- (7,3);
										\draw  [ultra thick, black]  (6,1) -- (6.5,1);
										\draw  [ultra thick, black]  (6.6,2.1) -- (7,2.1);
										\draw  [ultra thick, black]  (6.49,1) -- (6.49,2.5);
										\draw[ultra thick, black] (0,3) arc (90:270: 0.25 and 0.5);
										\draw[ultra thick, black] (7,3) arc (90:-90:0.25 and 0.435);
										\draw [black, very thick,->](0.1,1) -- (0.3,1);
										\fill[black] (0,1) circle (3pt);
										\fill[black] (6.49,2.5) circle (3pt);
										\draw[black, ultra thick, decorate, decoration={brace, amplitude=6pt, mirror}]
										(0.3,0.8) -- (5.8,0.8) node [midway, yshift=-0.4cm] {$2p-1$ crossings};
										
										\draw[ultra thick, red, dashed]  (1.7,1.5) circle[x radius=0.55, y radius=0.7];
										
									\end{scope}
								\end{tikzpicture}
							}
						\end{center}
						\caption{An $H(2)$-move transforming $\mathcal N(2p)$ into $U$.} \label{fig18}
					\end{figure}
					
					By Lemma~\ref{lemma3.8}, for an integer $p>1$, we have $\mathcal N(2p) \in [K]_p$, and from the proof of Lemma~\ref{lemma3.8} it is easy to see that $\mathcal N(2) \in [K]_{3/2}$. Thus, if $p>1$ is an integer or $p=3/2$, the knotoid class $[K]_p$ can be reduced to the class $[U]$ by a single $H(2)$-move. From Lemmas~\ref{lemma3.10} and~\ref{lemma3.14}, it follows that the knotoid $\mathcal L(m)$ can be transformed into $U$ by at most two $H(2)$-moves. By Lemma~\ref{lemma3.9}, for an integer $p>1$, we have $\mathcal L(p) \in [K]_{p-\frac{1}{2}}$. Therefore, if $p \in \frac{1}{2} \mathbb Z$ and $p \geq \frac{3}{2}$, the knotoid class $[K]_p$ can be transformed into $[U]$ by at most two $H(2)$-moves.
					
					Let $\delta_1$ be the minimum number of $H(2)$-moves required to transform $[K]_1$ into $[U]$, and let $\delta_{\frac{1}{2}}$ be the minimum number of $H(2)$-moves required to transform $[K]_{\frac{1}{2}}$ into $[U]$. Set $\delta = \max \{ 2, \delta_1, \delta_{\frac{1}{2}} \}$. Then the diameter of the graph $G(H(2), \mathbf{span}(V))$ does not exceed $2\delta$. Since $d_{H(2)} \geq d_{H(n)}$ for all $n\geq 3$, the diameter of the graph $G(H(n), \mathbf{span}(V))$ is also finite, and the graph is Gromov hyperbolic. The theorem is proved.
				\end{proof}
				
				\begin{remark} 
					As seen from the proof of Theorem~\ref{theorem1.6}, if the classes $[K]_1$ and $[K]_{\frac{1}{2}}$ were empty, one could assert that the diameter of the graph $G(H(n),\mathbf{span}(V))$ does not exceed $4$. However, to the best of our knowledge, the question of the existence of knotoids whose Jones polynomial span is equal to $1$ or $\frac{1}{2}$ remains open.
				\end{remark}
	\smallskip
	\noindent
	\textbf{Tool use.} The authors used AI for proofreading. The mathematical ideas and results presented in the article were conceives of and written by the authors without the aid of AI. The authors take full responsibility for all contents.

				\usetikzlibrary{decorations.markings}

				\end{document}